\documentclass[a4paper,11pt]{article}
\usepackage{srcltx}
\usepackage{indentfirst}
\usepackage{epsfig}
\usepackage{psfrag}
\usepackage{graphicx}
\usepackage{overpic}
\usepackage{multirow}
\usepackage{longtable}
\usepackage{hyperref}
\usepackage[hang]{subfigure}
\usepackage{amsmath,amssymb, amsthm}
\usepackage{color}
\usepackage[mathscr]{euscript}
\usepackage{cases}
\usepackage{enumerate}
\usepackage[toc,page,title,titletoc,header]{appendix} 
\usepackage{cite}
\usepackage{rotating}

\newcommand{\bg}{\boldsymbol{g}}
\newcommand{\bphi}{\boldsymbol{\phi}}
\newcommand{\bbphi}{\boldsymbol{\bar{\phi}}}
\newcommand{\bhphi}{\boldsymbol{\hat{\phi}}}

\newcommand\bI{\boldsymbol{I}}

\newcommand\bbR{\mathbb{R}}
\newcommand\bbN{\mathbb{N}}
\newcommand\bxi{\boldsymbol{\xi}}
\newcommand\bx{\boldsymbol{x}}

\newcommand\bq{\boldsymbol{q}}
\newcommand\bu{\boldsymbol{u}}
\newcommand\bv{\boldsymbol{v}}

\newcommand\bLambda{{\boldsymbol{\Lambda}}}

\newcommand\bF{\boldsymbol{F}}

\newcommand\dd{\,\mathrm{d}}
\newcommand\He{\mathit{He}}
\newcommand\Kn{\mathit{Kn}}

\newcommand\mH{\mathcal{H}}

\newcommand\mP{\mathcal{P}}
\newcommand\mQ{\mathcal{Q}}
\newcommand\bmP{\boldsymbol{\mP}}
\newcommand\bmQ{\boldsymbol{\mQ}}
\newcommand\bbmQ{\boldsymbol{\bar{\mQ}}}
\newcommand\bhmP{\boldsymbol{\hat{\mP}}}

\newcommand\mM{\mathcal{M}}

\newcommand\NRxx{NR$xx$~}

\newcommand\mF{\mathcal{F}}

\newcommand\sss{\scriptscriptstyle}

\graphicspath{{images/}}
{\theoremstyle{remark} \newtheorem{remark}{\bfseries Remark}
 \newtheorem{algorithm}{\bfseries Algorithm}
}

\title{Acceleration for Microflow Simulations of High-Order Moment
  Models by Using Lower-Order Model Correction}

\author{Zhicheng Hu \thanks{Department of Mathematics, College of
    Science, Nanjing University of Aeronautics and Astronautics,
    Nanjing 210016, China, email: {\tt huzhicheng@nuaa.edu.cn}},
  ~~ Ruo Li\thanks{HEDPS \& CAPT, LMAM \& School of Mathematical
    Sciences, Peking University, Beijing, China, email: {\tt
      rli@math.pku.edu.cn}.}, 
  ~~ Zhonghua Qiao \thanks{Department of Applied Mathematics, The Hong
    Kong Polytechnic University, Hung Hom, Kowloon, Hong Kong, email:
    {\tt zqiao@polyu.edu.hk}.}}

\begin{document}

\maketitle

\begin{abstract}
  We study the acceleration of steady-state computation for microflow,
  which is modeled by the high-order moment models derived recently
  from the steady-state Boltzmann equation with BGK-type collision
  term. By using the lower-order model correction, a novel nonlinear
  multi-level moment solver is developed. Numerical examples verify
  that the resulting solver improves the convergence significantly
  thus is able to accelerate the steady-state computation greatly. The
  behavior of the solver is also numerically investigated. It is shown
  that the convergence rate increases, indicating the solver would be
  more efficient, as the total levels increases. Three order reduction
  strategies of the solver are considered. Numerical results show that
  the most efficient order reduction strategy would be $m_{\sss l-1} =
  \lceil m_{\sss l} / 2 \rceil$.

\vspace*{4mm}
\end{abstract}

\section{Introduction} 
\label{sec:intro} 

Microflow simulations are of great interest in a number of high-tech
fields such as the Micro-Electro-Mechanical-Systems (MEMS) devices. As
the characteristic length shrinks into micro-scale regime, typically
ranging from $0.1\,\rm{\mu m}$ to several tens of microns, the
traditional Navier-Stokes-Fourier (NSF) model becomes frequently to
show large deviations from the real physics, and consequently one has
to find new models to simulate the microflows. Indeed, as the
fundamental equation of the kinetic theory, the Boltzmann equation is
able to describe flows well in such micro-scale regimes
\cite{Struchtrup}. Whereas, due to the intrinsic high dimensionality,
numerical solution of the Boltzmann equation still remains a real
challenge, even when its complicated integral collision operator (see
e.g. \cite{Cowling}) is replaced by some simplified relaxational
operators, such as the Bhatnagar-Gross-Krook (BGK) model \cite{BGK},
the ellipsoidal statistical BGK (ES-BGK) model \cite{Holway}, the
Shakhov model \cite{Shakhov}, etc. On the other hand, the Boltzmann
equation contains a detailed microscopic description of flows while in
practice we are mainly interested in the macroscopic quantities of
physical meaning, which can be extracted by taking moments from the
distribution function. Therefore, it still has a great demand nowadays
to develop appropriate macroscopic transport models, also referred to
as extended hydrodynamic models, which could give a satisfactory
description of flows with a considerable reduction of computational
effort. The moment method, originally introduced by Grad \cite{Grad},
was considered as one of most powerful approaches to this end.

Recently, in view of the importance of hyperbolicity for a well-posed
model, a globally hyperbolic moment method, following the Grad moment
method with an appropriate closure ansatz, was proposed in \cite{Fan,
  Fan_new}. Therein a series of high-order moment models, that are all
globally hyperbolic, is derived from the Boltzmann equation in a
systematic way. These models are viewed as extensions of the NSF model
in a macroscopic point of view, and the systematic derivation makes it
possible to use the model up to arbitrary order for practical
applications. From numerical point of view, they actually constitute a
semi-discretization of the Boltzmann equation, wherein the velocity
space is discretized by a certain Hermite spectral method. Benefit
from this, convergence of these models to the underlying Boltzmann
equation is expected with a high-order rate as the order of the model
increases, see \cite{Microflows1D} for example. Through a further
investigation of these hyperbolic moment models, a routine procedure
to derive globally hyperbolic moment models from general kinetic
equations was introduced in \cite{framework}.

To simulate flows by using the high-order moment models, an
accompanying numerical method, abbreviated as the \NRxx method, has
been developed in \cite{NRxx,NRxx_new,Cai,Li,Microflows1D}. It has a
uniform framework for the model of arbitrary order, thus the
implementation of the algorithm for the model of a large order would
not be encountered difficulties. Some successful applications not
limited to gas flow problems can be found in \cite{Wang, Hu2012}.
However, it turns out that the general designed \NRxx method becomes
inefficient, when the steady-state problems are considered or the
model of a sufficient large order is employed. While on the other
hand, there are quite some important applications in microflows, in
which the main concern is the steady-state solution, or the model of a
very large order is necessary for numerical purpose, see
e.g. \cite{Microflows1D}. In such situations, any improvement in
efficiency is worth of consideration.

As one of popular acceleration techniques for steady-state
computation, multigrid methods \cite{brandt2011book,
  hackbusch1985book} have been received increased attention in the
past few decades, and have been successfully applied in the classical
hydrodynamics \cite{li2008multigrid, hu2010robust,
  mavriplis2002assessment}. In our previous paper \cite{hu2014nmg}, 
a nonlinear multigrid (NMG) iteration, for the steady-state
computation of the hyperbolic moment models, has been developed, by
using the spatial coarse grid correction. Following the general design
idea of the \NRxx method, this NMG iteration is also unified for the
model of arbitrary order. It has been shown that significant
improvement in convergence has been obtained by the resulting NMG
solver in comparison to the direct time-stepping \NRxx scheme. Yet it
still takes a number of computational cost when the order of the model
is considerable large.

In this paper, we would consider the acceleration strategy for the
steady-state computation of the moment models from a novel direction.
It is pointed out that the hyperbolic moment models are in some sense
hierarchical models with respect to the model's order. Precisely
speaking, all equations in a moment model are contained in a
higher-order moment model, after removing the closure
ansatz. Observing this, it might be feasible to accelerate the
computation of the high-order moment model by using the lower-order
model correction, providing that the transformation operators between
moment models of different orders are appropriately proposed. This
would give rise to a multi-level moment algorithm for the high-order
moment model, as the NMG algorithm by using the spatial coarse grid
correction. The expectation, that such a new idea should be effective,
is mainly based on the following observations. Firstly, the
lower-order model correction can be viewed as the coarse grid
correction of velocity space, recalling that the moment model is
derived from the velocity discretization of the Boltzmann
equation. Consequently, the resulting multi-level moment solver would
constitute a multigrid solver of velocity space for the Boltzmann
equation. To the best of our knowledge, there is almost no efforts on
developing multigrid method of velocity space for the Boltzmann
equation in the literatures. Secondly, since a certain Hermite
spectral method is employed to derive the moment model from the
Boltzmann equation, the present multi-level moment solver would to
some extent coincide with the so-called $p$-multigrid method
\cite{fidkowski2005pmg, helenbrook2008solving} or spectral multigrid
method \cite{ronquist1987spectral, maday1988spectral}, which has been
successfully applied in various fields, see e.g. \cite{luo2008fast,
  kannan2011implicit, mascarenhas2010coupling-pmg, speck2015multi,
  wallraff2015multigrid}. Finally, numerical examples carried in the
present paper verify that this new idea is indeed able to accelerate
the steady-state computation significantly.

To accomplish the multi-level moment algorithm, the framework of
nonlinear multigrid algorithm developed in \cite{hackbusch1985book}
would be used. The implementation follows the basic idea of the \NRxx
method, such that the resulting nonlinear multi-level moment (NMLM)
solver also has a uniform framework for the model of arbitrary order,
and has the same input and output interfaces as the NMG solver
introduced in \cite{hu2014nmg}. Moreover, the transformation operators
between models of different orders could be implemented efficiently
under the framework of the \NRxx method. For the smoother of the NMLM
solver, the Richardson iteration with a cell-by-cell symmetric
Gauss-Seidel acceleration is proposed. A remaining important issue is
how to choose the order sequence for the NMLM solver, such that the
resulting solver not only improves the convergence rate but also saves
considerable computational cost. To this end, three order reduction
strategies are numerically investigated in the current paper to give
the best order reduction strategy. The behavior of the proposed NMLM
solver, with respect to the total levels of the solver, is also
numerically investigated. It turns out that the convergence rate is
improved as the total levels increases.

The remainder of this paper is organized as follows. The governing
Boltzmann equation and the corresponding hyperbolic moment models of
arbitrary order with a unified spatial discretization are briefly
reviewed in section \ref{sec:model}. Then the nonlinear multi-level
moment solver for the high-order moment model is comprehensively
introduced in section \ref{sec:method}. Its behavior is numerically
investigated in section \ref{sec:example} by two examples, which also
shows the robustness and efficiency of the proposed multi-level moment
solver. Finally, we give some concluding remarks in section
\ref{sec:conclusion}.


\section{The governing equations}
\label{sec:model}
In this section, we give a brief review of the governing Boltzmann
equation with BGK-type collision term in microflows, and the globally
hyperbolic moment models of arbitrary order, followed with a unified
spatial discretization.

\subsection{Boltzmann equation with BGK-type collision term}
\label{sec:model-boltzmann}
In the kinetic theory of microflows, the probability density of
finding a microscopic particle with velocity $\bxi\in \bbR^3$ at
position $\bx\in \Omega \subset \bbR^D$ $(D=1,2, \text{or } 3)$ is
measured by the distribution function $f(\bx,\bxi)$, whose evolution
is governed by the Boltzmann equation of the form
\begin{equation}
  \label{eq:boltzmann}
  \bxi \cdot \nabla_{\bx} f + \bF \cdot \nabla_{\bxi} f = Q(f)
\end{equation}
in the steady-state case. Here $\bF$ is the acceleration of particles
due to external forces, and the right-hand side $Q(f)$ is the
collision term representing the interaction between particles. As can
be seen in \cite{Cowling}, the original Boltzmann collision term is a
multi-dimensional integral, which turns out to be too inconvenient to
handle for numerical solution. Alternatively, several simplified
collision models are already able to capture the major physical
features of interest in a number of cases. In the present work, we
focus on the class of simplified relaxation models for $Q(f)$, saying,
the BGK-type collision term, which has a uniform expression given by
\begin{align}
  \label{eq:col-relaxation}
  Q(f) = \nu (f^{\text{E}} - f), 
\end{align}
where $\nu$ is the average collision frequency that is assumed
independent of the particle velocity, and $f^{\text{E}}$ is the
equilibrium distribution function depending on the specific model
selected. For instance, we have:
\begin{itemize}
\item For the ES-BGK model \cite{Holway}, $f^{\text{E}}$ is an
  anisotropic Gaussian distribution
  \begin{equation}
    \label{eq:ES-f^E}
    f^{\text{E}}(\bx,\bxi) = \frac{\rho(\bx)}{m_* \sqrt{ \det[2\pi \bLambda(\bx)]}
    } \exp\left(-\frac{1}{2} ( \bxi - \bu(\bx) )^T [\bLambda(\bx)]^{-1}
      (\bxi-\bu(\bx)) \right),
  \end{equation}
  where $\bLambda=(\lambda_{ij})$ is a $3\times 3$ matrix with
  \begin{align*}
    \lambda_{ij}(\bx) = \theta(\bx) \delta_{ij} + \left( 1-
      \frac{1}{\Pr} \right) \frac{\sigma_{ij}(\bx)}{\rho(\bx)}, \quad i,
    j = 1,2,3.
  \end{align*}
\item For the Shakhov model \cite{Shakhov}, $f^{\text{E}}$ reads 
  \begin{align}
    \label{eq:Shakhov-f^E}
    & f^{\text{E}}(\bx, \bxi) = \left[ 1+\frac{(1-\Pr)(\bxi-\bu(
        \bx))\cdot \bq( \bx) }{5\rho(\bx) [\theta(\bx)]^2} \left(
        \frac{ \vert \bxi-\bu(\bx) \vert^2}{\theta(\bx)} - 5
      \right)\right] f^{\text{M}}(
    \bx, \bxi),
  \end{align}
  where $f^M$ is the local Maxwellian given by 
  \begin{align}
    \label{eq:BGK-f^E}
    f^{\text{M}}(\bx,\bxi) = \frac{\rho( \bx)}{m_* [2\pi
      \theta(\bx)]^{3/2}} \exp\left(-\frac{\vert \bxi -
        \bu(\bx) \vert^2}{2\theta(\bx)} \right).
  \end{align}
\end{itemize}
In the above equations, $m_*$ is the mass of a single particle,
$\delta_{ij}$ is the Kronecker delta symbol, and $\Pr$ is the Prandtl
number. The macroscopic quantities, i.e., density $\rho$, mean
velocity $\bu$, temperature $\theta$, stress tensor $\sigma$, and heat
flux $\bq$, are related with the distribution function $f$ by
\begin{align}
  \label{eq:moments}
\begin{aligned}
  & \rho(\bx) = m_* \int_{\bbR^3} f(\bx,\bxi) \dd \bxi, \quad
  \rho(\bx) \bu(\bx) = m_* \int_{\bbR^3} \bxi f(\bx,\bxi) \dd \bxi, \\
  & \rho(\bx) \vert \bu(\bx) \vert^2 + 3 \rho(\bx) \theta(\bx) = m_*
  \int_{\bbR^3} \vert \bxi \vert^2 f(\bx, \bxi) \dd \bxi, \\ &
  \sigma_{ij}(\bx) = m_* \int_{\bbR^3} (\xi_i - u_i(\bx))(\xi_j -
  u_j(\bx)) f(\bx,\bxi) \dd \bxi - \rho(\bx) \theta(\bx) \delta_{ij},
  \quad i,j = 1,2,3,\\ & \bq(\bx) = \frac{m_*}{2} \int_{\bbR^3} \vert
  \bxi-\bu(\bx) \vert^2 (\bxi-\bu(\bx)) f(\bx,\bxi) \dd \bxi.
\end{aligned}
\end{align}
Note in the special case $\Pr=1$, both the ES-BGK model and the
Shakhov model reduce to the simplest BGK model \cite{BGK}, for which
$f^{\text{E}} \equiv f^{\text{M}}$.

\subsection{Hyperbolic moment equations of arbitrary order}
\label{sec:model-HME}

For convenience, we introduce $\mF^{[\tilde{\bu},\tilde{\theta}]}$ and
$\mF_M^{[\tilde{\bu},\tilde{\theta}]}$ to denote, respectively, the
linear spaces spanned by Hermite functions
\begin{align}
  \label{eq:base}
  \mH_{\alpha}^{[\tilde{\bu}, \tilde{\theta}]}(\bxi) =
  \frac{1}{m_*(2\pi \tilde{\theta})^{^{3/2}}
    \tilde{\theta}^{^{|\alpha|/2}}} \prod\limits_{d=1}^3
  \He_{\alpha_d}({v}_d)\exp \left(-{{v}_d^2}/{2} \right),
  \quad {\bv} = \frac{\bxi-\tilde{\bu}}{
    \sqrt{\tilde{\theta}}},~\forall \bxi\in\bbR^3,
\end{align}
for all $\alpha \in \bbN^3$ and for $\alpha$ with $|\alpha| \leq M$,
where $M\geq 2$ is a positive integer, $[\tilde{\bu}, \tilde{\theta}]\in
\bbR^3\times\bbR^+$ are two parameters, $|\alpha|$ is the sum of all
its components given by $|\alpha| = \alpha_1 + \alpha_2 + \alpha_3$,
and $\He_n(\cdot)$ is the Hermite polynomial of degree $n$, i.e.,
\begin{equation*}
  \He_n(x) = (-1)^n\exp \left( x^2/2 \right) \frac{\dd^n}{\dd x^n} 
  \exp \left(-x^2/2 \right).
\end{equation*}
It is easy to show that all $\mH_{\alpha}^{[\tilde{\bu},
  \tilde{\theta}]} (\bxi)$ are orthogonal to each other over $\bbR^3$
with respect to the weight function $\exp\left(|{\bv}|^2/2\right)$. It
follows that $\mF_M^{[\tilde{\bu}, \tilde{\theta}]}$ forms a finite
dimensional subspace of $L^2\left(\bbR^3,\exp \left(|{\bv}|^2 / 2
  \right)\dd \bxi\right)$ with $\mF_M^{[\tilde{\bu},\tilde{\theta}]}
\subset \mF_{M+1}^{[\tilde{\bu},\tilde{\theta}]} \subset \cdots
\subset \mF^{[\tilde{\bu},\tilde{\theta}]}$.

Following the derivation of the hyperbolic moment system of an
arbitrary order $M$ presented in \cite{Fan_new, Li, Microflows1D}, the
distribution function $f$ is approximated in
$\mF_M^{[\tilde{\bu},\tilde{\theta}]}$ with the parameters
$\tilde{\bu}$ and $\tilde{\theta}$ are exactly the local mean velocity
$\bu(\bx)$ and temperature $\theta(\bx)$ determined from $f$ itself
via \eqref{eq:moments}, that is,
\begin{align}
  \label{eq:truncated-dis}
  f(\bx, \bxi) \approx \sum_{\vert \alpha \vert \leq M} f_\alpha(\bx)
  \mH_{\alpha}^{[\bu(\bx), \theta(\bx)]} (\bxi), 
\end{align}
With such an approximation, we have from \eqref{eq:moments} the
following relations
\begin{align}\label{eq:moments-relation}
  \begin{aligned}
    & f_0 = \rho, \qquad f_{e_1} = f_{e_2} = f_{e_3} = 0, \qquad
    \sum_{d=1}^3 f_{2e_d} = 0, \\ 
    & \sigma_{ij} = (1+\delta_{ij}) f_{e_i+e_j}, \quad q_i = 2
    f_{3e_i} + \sum_{d=1}^3 f_{2e_d+e_i}, \qquad i,j=1,2,3,
\end{aligned}
\end{align}
where $e_1$, $e_2$, $e_3$ are introduced to denote the multi-indices
$(1,0,0)$, $(0,1,0)$, $(0,0,1)$, respectively.
 
By plugging \eqref{eq:truncated-dis} into the Boltzmann equation
\eqref{eq:boltzmann} with BGK-type collision term
\eqref{eq:col-relaxation}, matching the coefficients of the same basis
function, and applying the regularization proposed in \cite{Fan_new},
the hyperbolic moment system of order $M$ is then obtained as follows
\begin{equation}
  \label{eq:mnt-eqs}
  \begin{split}
    & \sum_{j=1}^D \Bigg[ \left( \theta \frac{\partial f_{\alpha -
          e_j}}{\partial x_j} + u_j \frac{\partial
        f_{\alpha}}{\partial x_j} + (1-\delta_{|\alpha|,M})(\alpha_j +
      1) \frac{\partial f_{\alpha+e_j}}{\partial x_j} \right) \\ &+
    \sum_{d=1}^3 \frac{\partial u_d}{\partial x_j} \left( \theta
      f_{\alpha-e_d-e_j} + u_j f_{\alpha-e_d} +
      (1-\delta_{|\alpha|,M}) (\alpha_j + 1) f_{\alpha-e_d+e_j}
    \right) \\ &+ \frac{1}{2} \frac{\partial \theta}{\partial x_j}
    \sum_{d=1}^3 \left( \theta f_{\alpha-2e_d-e_j} + u_j
      f_{\alpha-2e_d} + (1-\delta_{|\alpha|,M}) (\alpha_j + 1)
      f_{\alpha-2e_d+e_j} \right) \Bigg] \\ &= \sum_{d=1}^3 F_d
    f_{\alpha-e_d} + \nu (f^{\text{E}}_\alpha - f_\alpha), \qquad
    |\alpha| \leq M,
    \end{split}
\end{equation}
where $F_d$ is the $d$th component of the acceleration $\bF$, and
$f^{\text{E}}_\alpha$ are coefficients of the projection of
$f^{\text{E}}$ in the same function space $\mF_M^{[\bu,\theta]}$,
namely,
\begin{align}
  \label{eq:f^E-expansion}
  f^{\text{E}}(\bx, \bxi) \approx \sum_{|\alpha| \leq M}
  f^{\text{E}}_{\alpha}(\bx) \mH_{\alpha}^{[\bu(\bx),\theta(\bx)]}
  (\bxi).
\end{align}
As can be seen in \cite{Li} and \cite{Microflows1D}, the coefficients 
$f_{\alpha}^{\text{E}}$ can be analytically calculated for the Shakhov
model and the ES-BGK model.

The moment system \eqref{eq:mnt-eqs} is usually regarded as
macroscopic transport model in the kinetic theory, while from the
derivation point of view, it is actually a semi-discretization of the
Boltzmann equation, where the velocity space is discretized by a
certain Hermite spectral method. Consequently, the moment system
\eqref{eq:mnt-eqs} is expected to converge to the underlying Boltzmann
equation with a high-order rate as the system's order $M$ increases,
when the solution is smooth. Meanwhile, it allows us to return to the
Boltzmann equation to construct unified numerical solvers for the
moment system of arbitrary order. In turn, any solver developed for
the moment system can also be viewed as a solver for the Boltzmann
equation.

From \eqref{eq:mnt-eqs} we see that all moments, including the mean
velocity $\bu$, the temperature $\theta$ and the expansion
coefficients $f_{\alpha}$, are nonlinearly coupled with each
other. With additional relations given in \eqref{eq:moments-relation},
we have that the total number of independent moments in
\eqref{eq:mnt-eqs} is equal to the number of equations, which is clear
to be
\begin{align}
  \label{eq:total-moments-number}
\mM_M = {M+3 \choose 3},
\end{align}
e.g., $\mM_{10} = 286$ and $\mM_{26} = 3654$. It turns out that the
system might be very large when a high-order moment model is under
consideration, implying the computational cost would be still
considerable for a general designed numerical method. While on the
other hand, high-order moment model such as $M=10$ is commonly used in
microflow simulations, as can be seen in \cite{Microflows1D}, where we
can even see that the hyperbolic moment model up to the order $M=26$
is necessary for the planar Couette flow with $\Kn=1.199$.

\subsection{Spatial discretization}
\label{sec:model-dis}

In the rest of this paper, we restrict ourselves to one spatial
dimensional case. A unified finite volume discretization for the
moment model \eqref{eq:mnt-eqs} of arbitrary order can be obtained by
the so-called \NRxx method, which was first introduced in \cite{NRxx,
  Cai} and then developed in \cite{Li,NRxx_new,
  Microflows1D}. Specifically, we begin with the spatial finite volume
discretization of the Boltzmann equation \eqref{eq:boltzmann}, which
can be written in a general framework of the form
\begin{align}
  \label{eq:mnt-eq-dis}
  \frac{F(f_i(\bxi),f_{i+1}(\bxi)) - F(f_{i-1}(\bxi),
    f_i(\bxi))}{\Delta x_i} = G(f_i(\bxi)),
\end{align}
over the $i$th grid cell $[x_i, x_{i+1}]$, where $\{x_i\}_{i=0}^N$
constitute a mesh of the spatial domain $[0,L]$ with the length of the
$i$th cell to be $\Delta x_i = x_{i+1} - x_{i}$. Here $f_i(\bxi)$ is
the approximate distribution function on the $i$th cell,
$F(f_i,f_{i+1})$ is the numerical flux defined at $x_{i+1}$, the right
boundary of the $i$th cell, and the right-hand side $G(f_i)$
corresponds to the discretization of the acceleration and collision
terms of the Boltzmann equation \eqref{eq:boltzmann}. With the
assumption that $f_i(\bxi)$ belongs to a function space
$\mF_M^{[\tilde{\bu}_i, \tilde{\theta}_i]}$, i.e.,
\begin{align}
  \label{eq:dis-expansion-i}
  f_i(\bxi) = \sum_{|\alpha |\leq M} f_{i,\alpha}
  \mH_{\alpha}^{[\tilde{\bu}_i, \tilde{\theta}_i]}(\bxi),
\end{align}
all terms of \eqref{eq:mnt-eq-dis}, numerical fluxes $F(f_{i-1},f_i)$,
$F(f_i,f_{i+1})$ and the right-hand side $G(f_i)$, can be computed
and approximated as the functions in the same space
$\mF_M^{[\tilde{\bu}_i, \tilde{\theta}_i]}$, that is, they can be
expressed in terms of the same basis functions of $f_i(\bxi)$ as
follows,
\begin{align}
  \label{eq:flux-force-expansion}
  \begin{aligned}
    & F(f_{i-1},f_i) = \sum_{|\alpha | \leq M} F_{\alpha}(f_{i-1},f_i)
    \mH_{\alpha}^{[\tilde{\bu}_i, \tilde{\theta}_i]} (\bxi), \\ &
    F(f_i,f_{i+1}) = \sum_{|\alpha | \leq M} F_{\alpha}(f_{i},f_{i+1})
    \mH_{\alpha}^{[\tilde{\bu}_i, \tilde{\theta}_i]}(\bxi), \\ &
    G(f_i(\bxi)) = \sum_{|\alpha | \leq M} G_{i,\alpha}
    \mH_{\alpha}^{[\tilde{\bu}_i, \tilde{\theta}_i]}(\bxi).
  \end{aligned}
\end{align}
Substituting the above expansions into \eqref{eq:mnt-eq-dis} and
matching the coefficients of the same basis function
$\mH_{\alpha}^{[\tilde{\bu}_i, \tilde{\theta}_i]}(\bxi)$, we then get
a system which is a discretization of the hyperbolic moment system
\eqref{eq:mnt-eqs} on the $i$th cell, providing that the parameters
$\tilde{\bu}_i$, $\tilde{\theta}_i$ are mean velocity $\bu_i$ and
temperature $\theta_i$ of the $i$th cell, respectively, such that the
relation \eqref{eq:moments-relation} holds for $f_{i,\alpha}$, and the
numerical flux $F(f_i,f_{i+1})$ is designed specially to coincide with
the hyperbolicity of the moment system. Accordingly, the set of mean
velocity $\bu_i$, temperature $\theta_i$ and expansion coefficients
$f_{i,\alpha}$ forms a solution of the moment system on the $i$th
cell. In the rest of this paper, we would equivalently say the
corresponding distribution function $f_i(\bxi) \in
\mF_M^{[\bu_i,\theta_i]}$ is the solution of the moment system on the
$i$th cell for simplicity.

From the moment system \eqref{eq:mnt-eqs}, we can easily deduce that
$G_{i,\alpha} = \sum_{d=1}^3 F_{i,d} f_{i,\alpha-e_d} + \nu_i
(f_{i,\alpha}^{\text{E}} - f_{i,\alpha})$, whereas the computation of
the numerical fluxes $F(f_{i-1}, f_{i})$ and $F(f_i, f_{i+1})$,
subsequently the coefficients $F_{\alpha}(f_{i-1}, f_{i})$ and
$F_{\alpha}(f_i, f_{i+1})$, is not straightforward. In order to
written the numerical fluxes in the form given in
\eqref{eq:flux-force-expansion}, it usually requires a transformation
between $\mF_M^{[\bu_{i\pm 1}, \theta_{i\pm 1}]}$ and $\mF_M^{[\bu_i,
  \theta_i]}$, no matter which kind of numerical flux is chosen, since
the solution $f_{i\pm 1}(\bxi)\in \mF_M^{[\bu_{i\pm 1}, \theta_{i \pm
    1}]}$ are originally expressed in terms of different set of basis
functions. A fast transformation between two spaces,
$\mF_M^{[\bu,\theta]}$ and $\mF_M^{[\tilde{\bu}, \tilde{\theta}]}$,
which constitutes the core of the \NRxx method, has already been
provided in \cite{NRxx}. In the current paper, we would call such
transformation, whenever it is necessary, without explicitly pointing
out. Additionally, the numerical flux presented in \cite{Microflows1D}
is employed in our experiments for comparison.


\section{Numerical methods}
\label{sec:method}
This section is devoted to develop an efficient solver for the
high-order moment model \eqref{eq:mnt-eqs}, following the idea to
accelerate the computation by using the lower-order moment model
correction. We first introduce a basic iterative method for the moment
model \eqref{eq:mnt-eqs} of a given order upon the unified
discretization \eqref{eq:mnt-eq-dis}, then illustrate the key
ingredients of using the lower-order model correction, and finally
give a multi-level moment solver for the high-order moment model
\eqref{eq:mnt-eqs}.

\subsection{Basic nonlinear iteration}
\label{sec:method-sgs}
Defining the local residual on the $i$th cell by
\begin{align}
  \label{eq:residual}
  R_i(f_{i-1}, f_{i}, f_{i+1}) = \frac{F(f_i(\bxi),f_{i+1}(\bxi)) -
    F(f_{i-1}(\bxi), f_i(\bxi))}{\Delta x_i} - G(f_i(\bxi)),
\end{align}
the discretization \eqref{eq:mnt-eq-dis} can be rewritten into 
\begin{align}
  \label{eq:residual-eq-0}
  R_i(f_{i-1}, f_{i}, f_{i+1}) = r_i(\bxi),
\end{align}
with $r_i(\bxi) \equiv 0$, where $r_i(\bxi) \in \mF_M^{[\bu_i,
  \theta_i]}$ is a known function in a slightly more general sense. It
is apparent that the above discretization gives rise to a nonlinear
system coupling all unknowns, i.e., $\bu_i$, $\theta_i$ and
$f_{i,\alpha}$, $i=0,1,\ldots,N-1$, $|\alpha | \leq M$,
together. Since the discretization relies on the basis functions which
usually change on different cells, it is quite difficult to obtain a
global linearization for the discretization problem. Alternatively, we
consider a localization strategy of using the cell-by-cell
Gauss-Seidel method.

A symmetric Gauss-Seidel (SGS) iteration, to produce a new approximate
solution $f^{n+1}$ with $f_i^{n+1}(\bxi) \in \mF_M^{[\bu_i^{n+1},
  \theta_i^{n+1}]}$ from a given approximation $f^{n}$ with
$f_i^{n}(\bxi) \in \mF_M^{[\bu_i^{n}, \theta_i^{n}]}$, consists of two
loops in opposite directions as follows.
\begin{enumerate}
\item Loop $i$ increasingly from 0 to $N-1$, and obtain
  $f_i^{n+\frac{1}{2}}(\bxi)$ by solving
  \begin{align}
    \label{eq:gs-interval}
    R_i(f_{i-1}^{n+\frac{1}{2}},f_i^{n+\frac{1}{2}},f_{i+1}^n)= r_i(\bxi). 
  \end{align}
\item Loop $i$ decreasingly from $N-1$ to $0$, and obtain
  $f_i^{n+1}(\bxi)$ by solving
  \begin{align}
    \label{eq:gs-reverse-interval}
    R_i(f_{i-1}^{n+\frac{1}{2}},f_i^{n+1},f_{i+1}^{n+1})= r_i(\bxi). 
  \end{align}
\end{enumerate}
The Gauss-Seidel method reduces the original global problem into a
sequence of local problems, i.e., \eqref{eq:gs-interval} or
\eqref{eq:gs-reverse-interval}, on each cell with the distribution
function on that cell as the only unknown.  Thereby, both
\eqref{eq:gs-interval} and \eqref{eq:gs-reverse-interval} can be
abbreviated to
\begin{align}
  \label{eq:residual-eq}
  R_i(f_{i}) = r_i(\bxi),
\end{align}
by removing the superscripts and the dependence on the distribution
function $f_{i-1}(\bxi)$, $f_{i+1}(\bxi)$ on the adjacent
cells. Certainly, the equation \eqref{eq:residual-eq} is still a
nonlinear problem. In \cite{hu2014nmg}, a Newton type method has been
proposed to solve it, wherein numerical differentiation was employed
to calculate the Jacobian matrix instead of the complicated analytical
derivation. The resulting iteration, the so-called SGS-Newton
iteration, exhibits faster convergence rate than a common explicit
time-integration scheme. Through a number of numerical tests, however,
we observed that for a general code implementation, the computational
cost of each SGS-Newton iteration grows rapidly as the system's order
$M$ increases, leading that the total cost might be more expensive
than the explicit time-integration method for a sufficient high-order
moment model. Although optimization of the implementation of numerical
differentiation can improve the efficiency of the method, such an
optimization does usually heavily rely on the specific choice of the
numerical flux, hence loses the generality of the method.

Currently, we are focusing on establishing the framework and verifying
the effectiveness of the idea using the lower-order model correction
to accelerate the computation of the high-order moment model. So we
solve \eqref{eq:residual-eq} in this paper by one step of a simple
relaxation method, namely, Richardson iteration, as in
\cite{hu2015}. The Richardson iteration reads
\begin{align}
  \label{eq:richardson}
  f_i^{n+1}(\bxi) = f_i^{n}(\bxi) + \omega_i \left(r_i(\bxi) -
    R_i(f_{i}^{n}) \right), 
\end{align}
which numerically consists of two steps as follows:
\begin{enumerate}
\item Compute an intermediate distribution function $f_i^*(\bxi)$ in
  $\mF_M^{[\bu_i^{n},\theta_i^{n}]}$, that is, its expansion
  coefficients $f_{i,\alpha}^*$ in terms of the basis functions
  $\mH_{\alpha}^{[\bu_i^n, \theta_i^n]} (\bxi)$ are calculated by
  \begin{align*}
    f_{i,\alpha}^* = f_{i,\alpha}^{n} + \omega_i \left(r_{i,\alpha} -
      R_{i,\alpha}\right), \quad |\alpha | \leq M,
  \end{align*}
  where $f_{i,\alpha}^{n}$, $r_{i,\alpha}$, and $R_{i,\alpha}$
  represent expansion coefficients respectively of
  $f_{i,\alpha}^{n}(\bxi)$, $r_i(\bxi)$ and $R_i(f_i^{n})$ in terms of
  the same basis functions.
\item Compute the new macroscopic velocity $\bu_i^{n+1}$ and
  temperature $\theta_i^{n+1}$ from $f_i^*(\bxi)$, then project
  $f_i^*(\bxi)$ into $\mF_M^{[\bu_i^{n+1},\theta_i^{n+1}]}$ to obtain
  $f_i^{n+1}(\bxi)$.
\end{enumerate}
The relaxation parameter $\omega_i$ in \eqref{eq:richardson} is
selected according to the local CFL condition
\begin{align}
  \label{eq:local-CFL}
  \omega_i \frac{\lambda_{\max,i}}{\Delta x_i} < 1,
\end{align}
and the strategy to preserve the positivity of the local density and
temperature, see \cite{hu2014nmg} for details. Here,
$\lambda_{\max,i}$ is the largest value among the absolute values of
all eigenvalues of the hyperbolic moment model \eqref{eq:mnt-eqs} on
the $i$th cell.

Now we have a basic nonlinear iteration, referred to as SGS-Richardson
iteration in the rest of this paper, for the moment model
\eqref{eq:mnt-eqs} of a certain order. A single level solver would
then be obtained by performing this basic iteration until the steady
state has been achieved. The criterion indicating the steady state is
adopted as
\begin{align}
  \label{eq:steady-criterion}
  \left\Vert \tilde{R} \right\Vert \leq \mathit{tol},
\end{align}
where $\mathit{tol}$ is a given tolerance, and $\left \Vert \tilde{R}
\right \Vert$ is the norm of the global residual $\tilde{R}$ given by
\begin{align}
  \label{eq:global-l2-norm}
  \left\Vert \tilde{R} \right \Vert = \sqrt{ \frac{1}{L}\left(
      \sum_{i=0}^{N-1} \left\Vert \tilde{R}_{i} \right\Vert^2
      \Delta x_i \right)}. 
\end{align}
Here, the local residual $\tilde{R}_i$ is defined on the $i$th cell by
$\tilde{R}_i(\bxi) = r_i(\bxi) - R_i(f_{i-1},f_{i}, f_{i+1})$, and its
norm is computed using the weight $L^2$ norm of the linear space
$\mF_M^{[\bu_i,\theta_i]}$, that is,
\begin{align}
  \label{eq:l2-norm-dis-def}
  \left\Vert \tilde{R}_i \right \Vert = \sqrt{\int \left(
      \tilde{R}_i(\bxi) \right)^2 \exp\left(
      \frac{\left|\bxi-\bu_i \right|^2}{2\theta_i} \right)
    \dd \bxi}.
\end{align}
Using the orthogonality of basis functions, it follows that
\begin{align}
  \label{eq:l2-norm-dis}
  \left\Vert \tilde{R}_i \right \Vert = \sqrt{ \sum_{|\alpha | \leq M}
    C_\alpha \left| \tilde{R}_{i,\alpha} \right|^2}, 
\end{align}
where $C_\alpha = m_*^{-2}(2\pi)^{-3/2} \left( \theta_i
\right)^{-|\alpha|-3/2} \alpha!$ with $\alpha! = \alpha_1!  \alpha_2!
\alpha_3!$, and $\tilde{R}_{i,\alpha}$ is the expansion coefficients
of $\tilde{R}_i(\bxi)$ in $\mF_M^{[\bu_i,\theta_i]}$.
\begin{remark}
  It is not suitable to calculate \eqref{eq:l2-norm-dis} more simple
  with $C_\alpha = 1$, by noting that $f_{i,\alpha}$ as well as
  $\tilde{R}_{i,\alpha}$ has the same dimension unit with $\rho_i
  \theta_i^{|\alpha|/2}$. In fact, $C_\alpha= m_*^{-2}(2\pi)^{-3/2}
  \left( \theta_i \right)^{-|\alpha|-3/2} \alpha!$ is also used to
  make each term in the summation of \eqref{eq:l2-norm-dis} have the
  same dimension unit $m_*^{-2}\rho_i^2\theta_i^{-3/2}$. Perhaps it is
  better to replace the weight function $\exp( \left|\bxi-\bu_i
  \right|^2 / (2\theta_i) )$ in \eqref{eq:l2-norm-dis-def} by
  $\rho_i/(m_* f_i^M) = (2\pi\theta_i)^{3/2} \exp( \left|\bxi-\bu_i
  \right|^2 / (2\theta_i) )$, in the sense that now each term in the
  summation of \eqref{eq:l2-norm-dis} would be dimensionalized to
  $m_*^{-2}\rho_i^2$.
\end{remark}
\begin{remark}
  Limited by machine float-point precision, the calculation of
  \eqref{eq:l2-norm-dis} becomes inaccurate when $M$ is a little big,
  $M\geq 10$ for example. This influences the study on the performance
  of the proposed method in this paper. Noting on the other hand that
  the macroscopic quantities of physical interest can be obtained from
  the first several moments, we approximate the norm of the local
  residual by
  \begin{align}
  \label{eq:l2-norm-dis-p1}
  \left\Vert \tilde{R}_i \right \Vert \approx \sqrt{ \sum_{|\alpha |
      \leq \min\{M,3\}} C_\alpha \left| \tilde{R}_{i,\alpha}
    \right|^2},
\end{align}
instead of \eqref{eq:l2-norm-dis} in our numerical experiments.
The local residual computed by \eqref{eq:l2-norm-dis-p1}
  changes with respect to $M$ even for the same $\bu_i$, $\theta_i$,
  $f_{i,\alpha}$, $i=0,1,\ldots,N-1$, $|\alpha | \leq M$, and $M>4$,
  since the numerical flux presented in \cite{Microflows1D} depends on
  the eigenvalues of the moment model, which clearly change with
  respect to $M$. Therefore, we can still use
  \eqref{eq:l2-norm-dis-p1} to measure the residual and find the
  correct result for a high-order moment model, even when the
  steady-state solution of a lower-order moment model is employed as
  the initial value.  
\end{remark}

As explained in \cite{hu2015}, the SGS-Richardson iteration can be
viewed as the variation of an explicit time-integration scheme.
Consequently, although the total computational cost is saved a lot by
the SGS-Richardson iteration for it converges in general several times
faster than the explicit time-integration scheme, the asymptotic
behavior of both two methods are similar. For example, the increase
rate of the total iterations with respect to spatial grid number $N$
or model's order $M$ is similar for both the SGS-Richardson iteration
and the explicit time-integration scheme. In order to get a more
efficient solver, we have considered in \cite{hu2014nmg} and
\cite{hu2015} the strategy using the coarse grid correction to
accelerate the convergence, and it has been validated that the
resulting nonlinear multigrid solvers have a significant improvement
in efficiency.

In this paper, we would consider the acceleration strategy for the
high-order moment model from another direction. Precisely speaking, we
would like to accelerate the convergence by using the lower-order
model correction. The details for this new strategy will given in
the following subsections.

\subsection{Lower-order model correction}
\label{sec:method-lower}

Let us rewrite the underlying problem resulting from
\eqref{eq:residual-eq-0} of a high order $M$ into a global form as
\begin{align}
  \label{eq:high-order-problem}
  R_M(f_M) = r_{\sss M},
\end{align}
and suppose $\bar{f}_M$ with its $i$th component $\bar{f}_{M,i}(\bxi)
\in \mF_M^{[\bar{\bu}_{M, i}, \bar{\theta}_{M,i}]}$ is an approximate
solution for the above problem. Like with the spatial coarse grid
correction used in \cite{hu2014nmg}, the lower-order problem is given
by
\begin{align}
  \label{eq:lower-order-problem}
  R_m(f_m) = r_m \triangleq R_m(\tilde{I}_M^m \bar{f}_M) + I_M^m \left(r_{\sss M} -
    R_M(\bar{f}_M)\right),
\end{align}
where $\square_M^m$ is the restriction operators moving functions from
the high $M$th-order function space to a lower $m$th-order function
space. The lower-order operator $R_m$ is analogous to the high-order
counterpart $R_M$, that is, $R_m(f_m)$ is obtained by the
discretization formulation \eqref{eq:residual} of the $m$th-order
moment model. It follows that the lower-order problem
\eqref{eq:lower-order-problem} can be solved using the same strategy
as the high-order problem \eqref{eq:high-order-problem}. When the
solution $f_m$ of the lower-order problem
\eqref{eq:lower-order-problem} is obtained, the solution of the
high-order problem \eqref{eq:high-order-problem} is then corrected by
\begin{align}
  \label{eq:update-high-order}
  \hat{f}_M = \bar{f}_M + I_m^M\left( f_m - \tilde{I}_M^m \bar{f}_M
  \right),
\end{align}
where $I_m^M$ is the prolongation operator moving functions from the
$m$th-order function space to the $M$th-order function space.

Recalling that the moment model \eqref{eq:mnt-eqs} is derived from the
Boltzmann equation \eqref{eq:boltzmann} by a special Hermite spectral
discretization of the velocity space, we conclude that the above
lower-order model correction is in fact a coarse grid correction of
velocity space. Furthermore, the idea using lower-order model
correction does to some extent coincide with the so-called
$p$-multigrid method \cite{fidkowski2005pmg, helenbrook2008solving},
which accordingly provides us with a reference to design the solver
for our purpose.

\subsection{Restriction and prolongation}
\label{sec:method-restriction}
 
In the current work, the lower-order problem
\eqref{eq:lower-order-problem} is defined on the same spatial mesh as
the high-order problem \eqref{eq:high-order-problem}. Therefore, it is
enough to give the definition of the restriction and prolongation
operators on an individual element of the spatial mesh. For
simplicity, the index $i$ of the spatial element is omitted in this
subsection.

By means of the unified expression \eqref{eq:dis-expansion-i} which
deals with all moments of the model as a whole, we can design the
restriction and prolongation operators following the $p$-multigrid
method \cite{fidkowski2005pmg, helenbrook2008solving}. Let $\bphi^M$
and $\bphi^m$ denote the column vectors of basis functions spanning
the $M$th-order space $\mF_M^{[\bu_{M}, \theta_{M}]}$ and the
$m$th-order space $\mF_m^{[\bu_{m}, \theta_{m}]}$, respectively. The
weighted $L^2$ projection of $\bphi^m$ in $\mF_M^{[\bu_{M},
  \theta_{M}]}$ is then given by $\bmP^T \bphi^M$, where $\bmP$ is a
$\mM_M \times \mM_m$ matrix defined as 
\begin{align}
  \label{eq:prolongation-matrix}
  \bmP = \left( \int \bphi^M \left(\bphi^M\right)^T \exp\left(
      \frac{|\bxi - \bu_{\sss M}|^2}{2\theta_{M}}\right) \dd \bxi
  \right)^{-1} \int \bphi^M (\bphi^m)^T \exp\left( \frac{|\bxi -
      \bu_{\sss M}|^2}{2\theta_{M}}\right) \dd \bxi.
\end{align}
Similarly, the weighted $L^2$ projection of $\bphi^M$ in
$\mF_m^{[\bu_{m}, \theta_{m}]}$ is given by $\bmQ^T \bphi^m$, where
$\bmQ$ is a $\mM_m \times \mM_M$ matrix defined as 
\begin{align}
  \label{eq:restriction-matrix}
  \bmQ = \left( \int \bphi^m \left(\bphi^m\right)^T \exp\left(
      \frac{|\bxi - \bu_{m}|^2}{2\theta_{m}}\right) \dd \bxi
  \right)^{-1} \int \bphi^m \left(\bphi^M\right)^T \exp\left( \frac{|\bxi -
      \bu_{m}|^2}{2\theta_{m}}\right) \dd \bxi.
\end{align}
Thus, the prolongation operator $I_m^M$ and the residual operator
$I_M^m$ can be defined, respectively, by the matrix $\bmP$ and its
transpose $\bmP^T$. That is, for the functions $g_m = \left(\bphi^m
\right)^T\bg_{m} \in \mF_m^{[\bu_{m}, \theta_{m}]}$ and $g_{\sss M} =
\left( \bphi^M \right)^T \bg_{\sss M} \in \mF_M^{[\bu_{M},
  \theta_{M}]}$, where the bold symbols $\bg_{\sss M}$ and $\bg_m$ are
the column vectors of the corresponding expansion coefficients
$g_{\sss M,\alpha}$ and $g_{m,\alpha}$, we have 
\begin{align*}
  I_m^M g_m = \left(\bphi^M \right)^T \bmP\bg_m, \qquad I_M^m g_{\sss
    M} = \left( \bphi^m \right)^T \bmP^T \bg_{\sss M}.
\end{align*}
Usually, the solution restriction operator $\tilde{I}_M^m$ does not
have to be the same as the residual restriction operator $I_M^m$,
and can be defined as $\tilde{I}_M^m g_{\sss M} = \left( \bphi^m
\right)^T \bmQ \bg_{\sss M}$.

In contrast to the $p$-multigrid method, unfortunately, the
computation of the matrices $\bmP$ and $\bmQ$ would be very expensive,
since $\bu_m$, $\theta_m$ are commonly not equal to $\bu_{\sss M}$,
$\theta_M$, and even all these values, consequently the basis
functions $\bphi^m$ and $\bphi^M$, have been changing throughout the
iterative procedure. Not only that, the exact matrices $\bmP$ and
$\bmQ$ are in fact unknown when the restriction operators
$\tilde{I}_M^m$ and $I_M^m$ are applied in
\eqref{eq:lower-order-problem}, for $\bu_m$ and $\theta_m$ can not be
obtained until \eqref{eq:lower-order-problem} has been solved.

To find the way out, let us return to the lower-order problem
\eqref{eq:lower-order-problem}. As stated in previous, all terms of
\eqref{eq:lower-order-problem}, in the initial discretization of each
element, are represented in terms of $\bbphi^m$, the basis functions
of $\mF_m^{[\bar{\bu}_{m}, \bar{\theta}_{m}]}$ that is determined by
the initial guess $\bar{f}_m$. Without any other information, a good
choice for $\bar{f}_m$ might be that it takes conservative quantities
the same as the high-order solution $\bar{f}_M \in
\mF_M^{[\bar{\bu}_M, \bar{\theta}_M]}$, that is,
\begin{align}\label{eq:principle-lower-solution}
  \int \bar{f}_m \varphi \dd \bxi = \int \bar{f}_M \varphi \dd \bxi,
  \quad \varphi = \left(1,~ \bxi,~\frac{1}{2}|\bxi|^2\right)^T.
\end{align}
It follows that $\bar{\bu}_m = \bar{\bu}_{\sss M}$ and $\bar{\theta}_m
= \bar{\theta}_M$, which indicate that $\bbphi^m$ coincides with the
first $\mM_m$ functions of $\bbphi^M$, the basis functions of
$\mF_M^{[\bar{\bu}_M, \bar{\theta}_M]}$. Using the orthogonality of
the basis functions, the special matrix $\bbmQ$, defined as
\eqref{eq:restriction-matrix} for $\bbphi^M$ and $\bbphi^m$, becomes
$\bbmQ = [\bI, \boldsymbol{0}]$, where $\bI$ is the identity matrix of
order $\mM_m$ and $\boldsymbol{0}$ represents the $\mM_m \times (\mM_M
- \mM_m)$ zero matrix. Noting that the initial guess $\bar{f}_m$ in
practice is always taken by $\tilde{I}_M^m \bar{f}_M$, we now define
the restriction operator $\tilde{I}_M^m$ from
$\mF_M^{[\bar{\bu}_M, \bar{\theta}_M]}$ into $\mF_m^{[\bar{\bu}_{m},
  \bar{\theta}_{m}]}$ as $\tilde{I}_M^m g_{\sss M} = \left( \bbphi^m
\right)^T \bbmQ \bg_{\sss M}$, that is, $\tilde{I}_M^m$ is just a
simple truncation operator that directly gets rid of the part in terms
of the basis functions $\mH_{\alpha}^{[\bar{\bu}_M,
  \bar{\theta}_M]}(\bxi)$ with $|\alpha| > m$. Since the high-order
residual is finally projected into $\mF_m^{[\bar{\bu}_{m},
  \bar{\theta}_{m}]}$ in \eqref{eq:lower-order-problem}, we define the
residual restriction operator $I_M^m$ the same as $\tilde{I}_M^m$.

When the correction step \eqref{eq:update-high-order} is performed, we
can first calculate the new velocity $\hat{\bu}_M$ and temperature
$\hat{\theta}_M$, then the prolongation operator $I_m^M$
from $\mF_m^{[\bu_{m}, \theta_{m}]}$ into $\mF_M^{[\hat{\bu}_M,
  \hat{\theta}_M]}$ can be applied as $I_m^M g_m = \left(\bhphi^M
\right)^T \bhmP\bg_m$, where $\bhphi^M$ is the basis functions of the
updated high-order solution space, and $\bhmP$ is the matrix defined
as \eqref{eq:prolongation-matrix} for $\bhphi^M$ and $\bphi^m$. To
implement the prolongation procedure efficiently, the lower-order
correction in $\mF_m^{[\bu_{m}, \theta_{m}]}$ is first retruncated
into $\mF_M^{[\bu_{m}, \theta_{m}]}$, then projected into
$\mF_M^{[\hat{\bu}_M, \hat{\theta}_M]}$ by the transformation proposed
in \cite{NRxx}. In other words, $\bhmP$ is computed by $\bhmP= \bmP_0
\bbmQ^T$ instead of direct computation by the formula
\eqref{eq:prolongation-matrix}, where $\bmP_0$ is the matrix
representation of the transformation between two spaces with the same
order.

\subsection{Multi-level moment solver}
\label{sec:method-velocity}

Obviously, the lower-order problem \eqref{eq:lower-order-problem}
itself can also be solved by the two-level method using a much
lower-order model correction. Recursively applying this two-level
strategy then gives rise to a nonlinear multi-level moment (NMLM)
iteration.

Let $m_{\sss l}$, $l=0,1,\ldots,L$, denote the order of the
$l$th-level problem, and satisfy $2\leq m_0 < m_1 <\cdots < m_{\sss
  L}$. Then the $(l+1)$-level NMLM iteration, denoted by
$f_{m_l}^{n+1} = \text{NMLM}_l(f_{m_l}^n, r_{m_l})$, is given in the
following algorithm.
\begin{algorithm} [Nonlinear multi-level moment (NMLM) iteration]~
\label{alg:nmg}
  \begin{enumerate}
  \item If $l=0$, call the lowest-order solver, which will be given
    later, to have a solution $f_{m_0}^{n+1}$; otherwise, go to the
    next step.
  \item Pre-smoothing: perform $s_1$ steps of the
    \text{SGS-Richardson} iteration beginning with the initial
    approximation $f_{m_l}^n$ to obtain a new approximation
    $\bar{f}_{m_l}$.
  \item Lower-order model correction:
    \begin{enumerate}
    \item Compute the high-order residual as $\bar{R}_{m_l} = r_{m_l}
      - R_{m_l}(\bar{f}_{m_l})$.
    \item Prepare the initial guess of the lower-order problem by the
      restriction operator $\tilde{I}_{m_l}^{m_{l-1}}$ as
      $\bar{f}_{m_{l-1}} = \tilde{I}_{m_l}^{m_{l-1}} \bar{f}_{m_l}$.
    \item Calculate the right-hand side of the lower-order problem
      \eqref{eq:lower-order-problem} as $r_{m_{l-1}} =
      I_{m_l}^{m_{l-1}} \bar{R}_{m_l} + R_{m_{l-1}}(
      \bar{f}_{m_{l-1}})$.
    \item Recursively call the NMLM algorithm (repeat $\gamma$ times
      with $\gamma=1$ for a so-called $V$-cycle, $\gamma=2$ for a
      $W$-cycle, and so on) as
      \begin{align*}
        \tilde{f}_{m_{l-1}} = \text{NMLM}_{l-1}^\gamma (\bar{f}_{m_{l-1}},
        r_{m_{l-1}}).
      \end{align*}
    \item Correct the high-order solution by $\hat{f}_{m_{l}} =
      \bar{f}_{m_l}+ I_{m_{l-1}}^{m_{l}}(\tilde{f}_{m_{l-1}} - \bar{f}_{m_{l-1}})$.
    \end{enumerate}
  \item Post-smoothing: perform $s_2$ steps of the
    \text{SGS-Richardson} iteration beginning with $\hat{f}_{m_{l}}$
    to obtain the new approximation $f_{m_l}^{n+1}$.
  \end{enumerate}
\end{algorithm}
The $(l+1)$-level NMLM solver for the problem of order $m_{\sss l}$ is
then obtained by performing the above $(l+1)$-level NMLM iteration
until the steady state has been achieved. Obviously, the one-level
NMLM solver is just the single level solver of SGS-Richardson iteration.

Since the lowest-order problem is still a nonlinear problem with the
lowest-order operator $R_{m_0}$ analogous to the operator $R_{m_l}$ on
other order levels, a direct method for its exact solution is clearly
unavailable, and the SGS-Richardson iteration using as the smoothing
operator is again applied to give the lowest-order solver. In view of
that the spatial mesh is unchanged in the above NMLM algorithm,
accurately solving the lowest-order problem would lead to too much
SGS-Richardson iterations to make the whole NMLM solver
inefficient. Hence, only $s_3$ steps of the SGS-Richardson iteration
is performed in each calling of the lowest-order solver, where $s_3$
is a positive integer a little larger than the smoothing steps
$s_1+s_2$.

A remaining technical issue is how to set the order of the lower-order
problem. The order reduction strategy of either $m_{\sss l-1} =
m_{\sss l}-1$ or $m_{\sss l-1} = \lceil m_{\sss l} / 2 \rceil$ is
frequently used in the $p$-multigrid algorithm. Apart from them, the
strategy of $m_{\sss l-1} = m_{\sss l} - 2$ is also considered by
noting that the solution in our experiments exhibit a property
depending on the parity of the order of the model. In next section, we
will investigate the performance of all these three order reduction
strategies, and try to give the best one in the interest of improving
efficiency.


\section{Numerical examples}
\label{sec:example}
We present in this section two numerical examples, the planar Couette
flow and the force driven Poiseuille flow, to investigate the main
features of the proposed NMLM solver. For simplicity, we consider the
dimensionless case and the particle mass $m_*$ is always $1$. A
$V$-cycle NMLM solver with $s_1=s_2=2$ and $s_3=10$ is performed for
all numerical tests. The tolerance indicating the achievement of
steady state is set as $\mathit{tol}=10^{-8}$. We have observed that
the behavior of the NMLM solver are similar for the BGK-type collision
models. Thus only results for the ES-BGK collision model with the
Prandtl number $\Pr=2/3$ are given below.

To complete the problem, the Maxwell boundary conditions derived in
\cite{Li} are adopted for our moment models. As mentioned in
\cite{hu2014nmg}, such boundary conditions could not determine a
unique solution for the steady-state moment model
\eqref{eq:mnt-eqs}. In order to recover the consistent steady-state
solution with the time-stepping scheme and the NMG solver proposed in
\cite{hu2014nmg}, the correction employed in \cite{Mieussens2004,
  hu2014nmg} is also applied to the solution at each NMLM iterative
step.

\subsection{The planar Couette flow}
\label{sec:num-ex-couette}

The planar Couette flow is frequently used as benchmark test in the
microflows. Consider the gas in the space between two infinite
parallel plates, which have the same temperature $\theta^W$ and are
separated by a distance $L$. One plate is stationary, and the other is
translating with a constant velocity $u^W$ in its own plane. Although
there is no external force acting on the gas, that is, $\bF\equiv 0$,
the gas will still be driven by the motion of the plate, and finally
reach a steady state.

\begin{figure}[!htb]
{  \centering
  \subfigure[Density, $\rho$]{\includegraphics[width=0.231\textwidth]{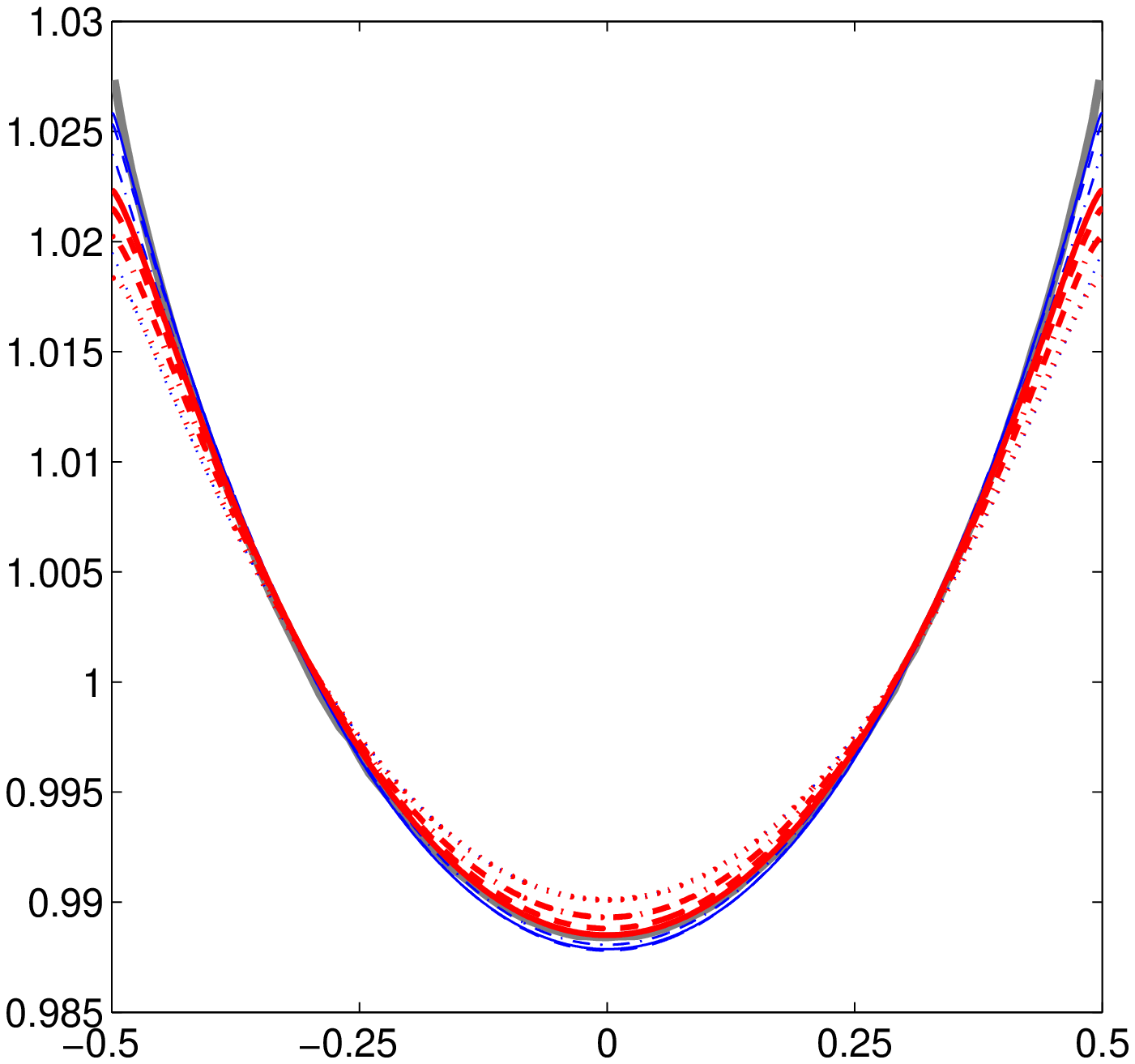}}
  \subfigure[Temperature, $\theta$]{\includegraphics[width=0.231\textwidth]{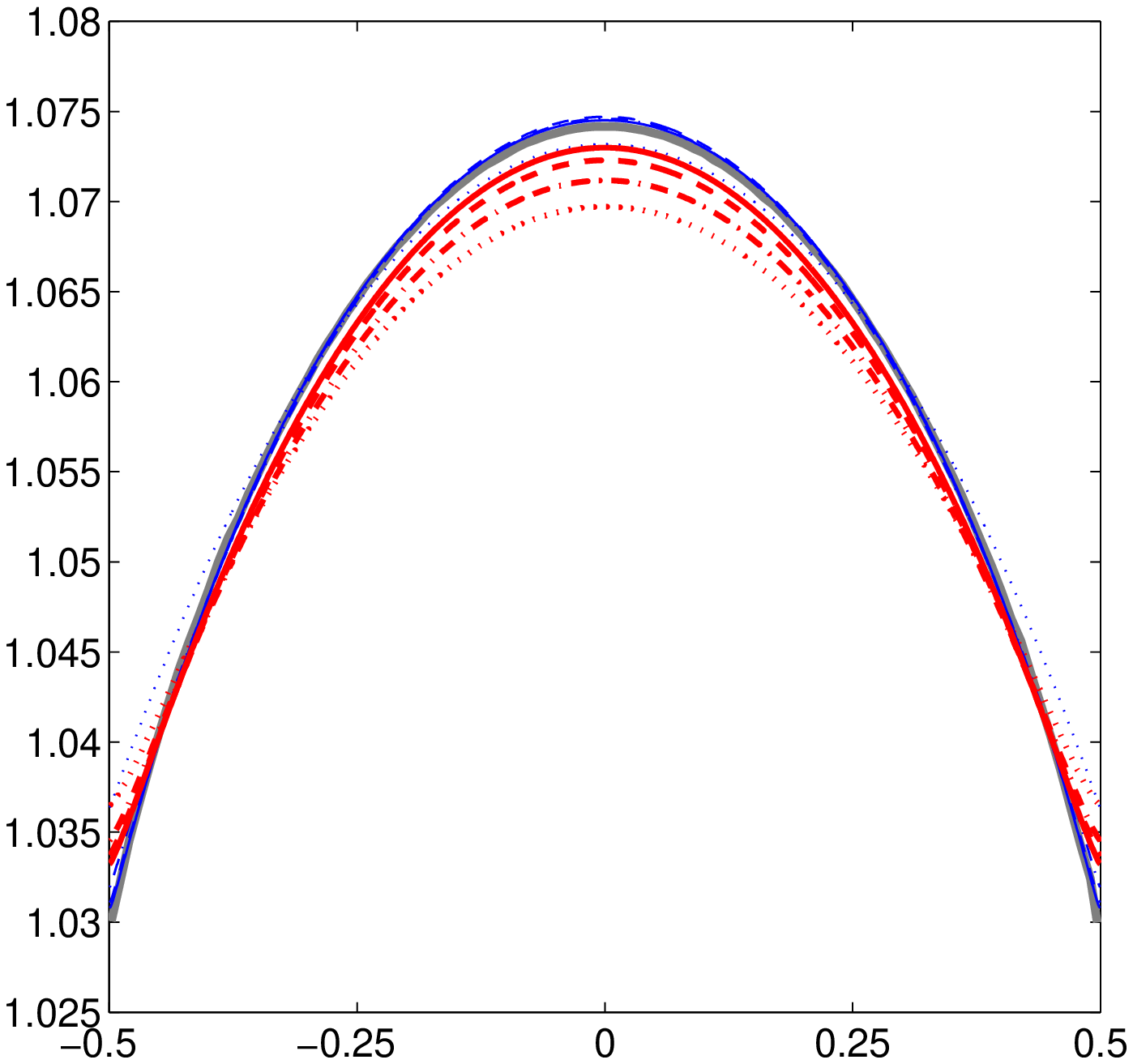}} 
  \subfigure[Shear stress, $\sigma_{12}$]{\includegraphics[width=0.236\textwidth]{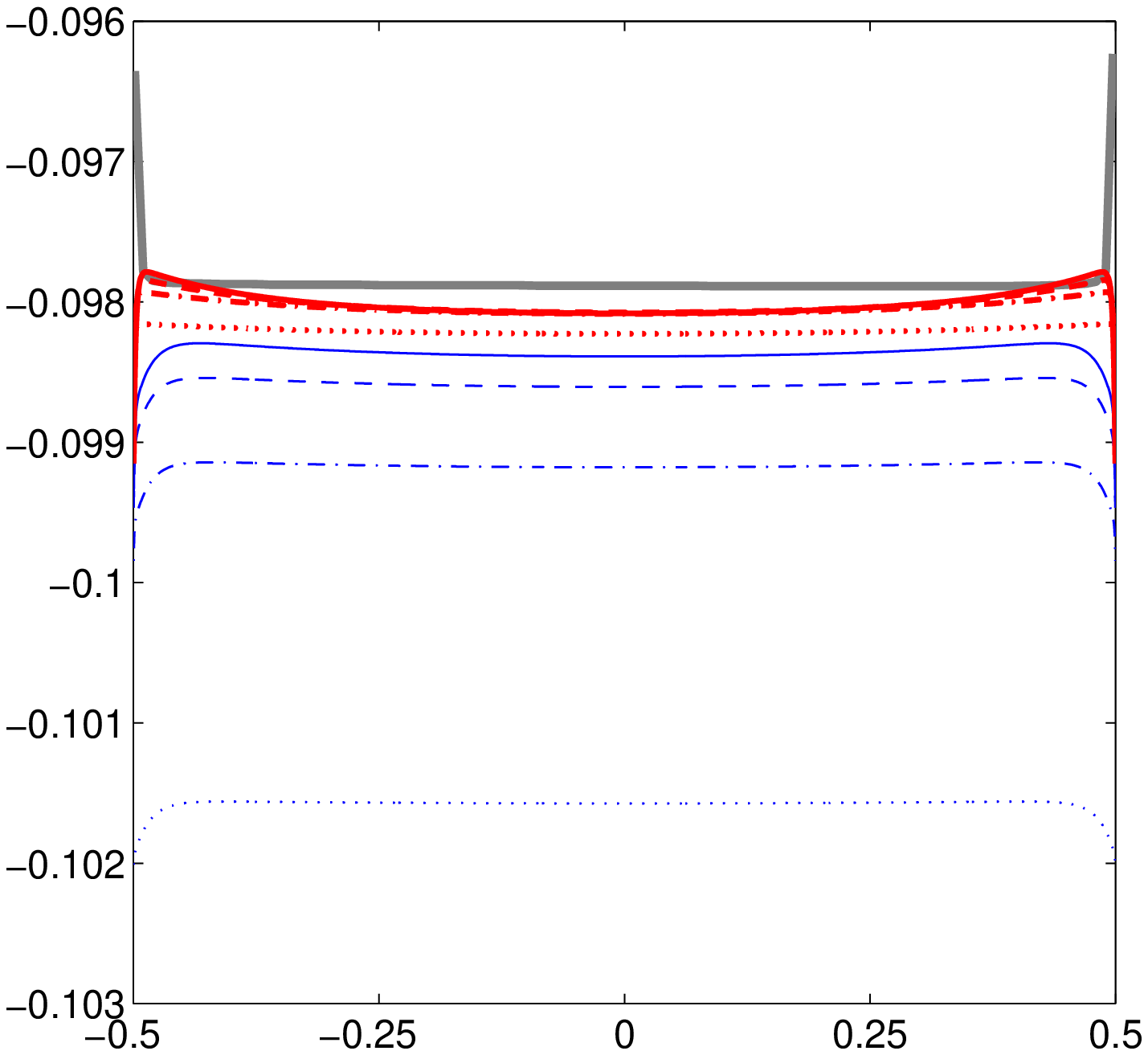}} 
  \subfigure[Heat flux, $q_1$]{\includegraphics[width=0.278\textwidth]{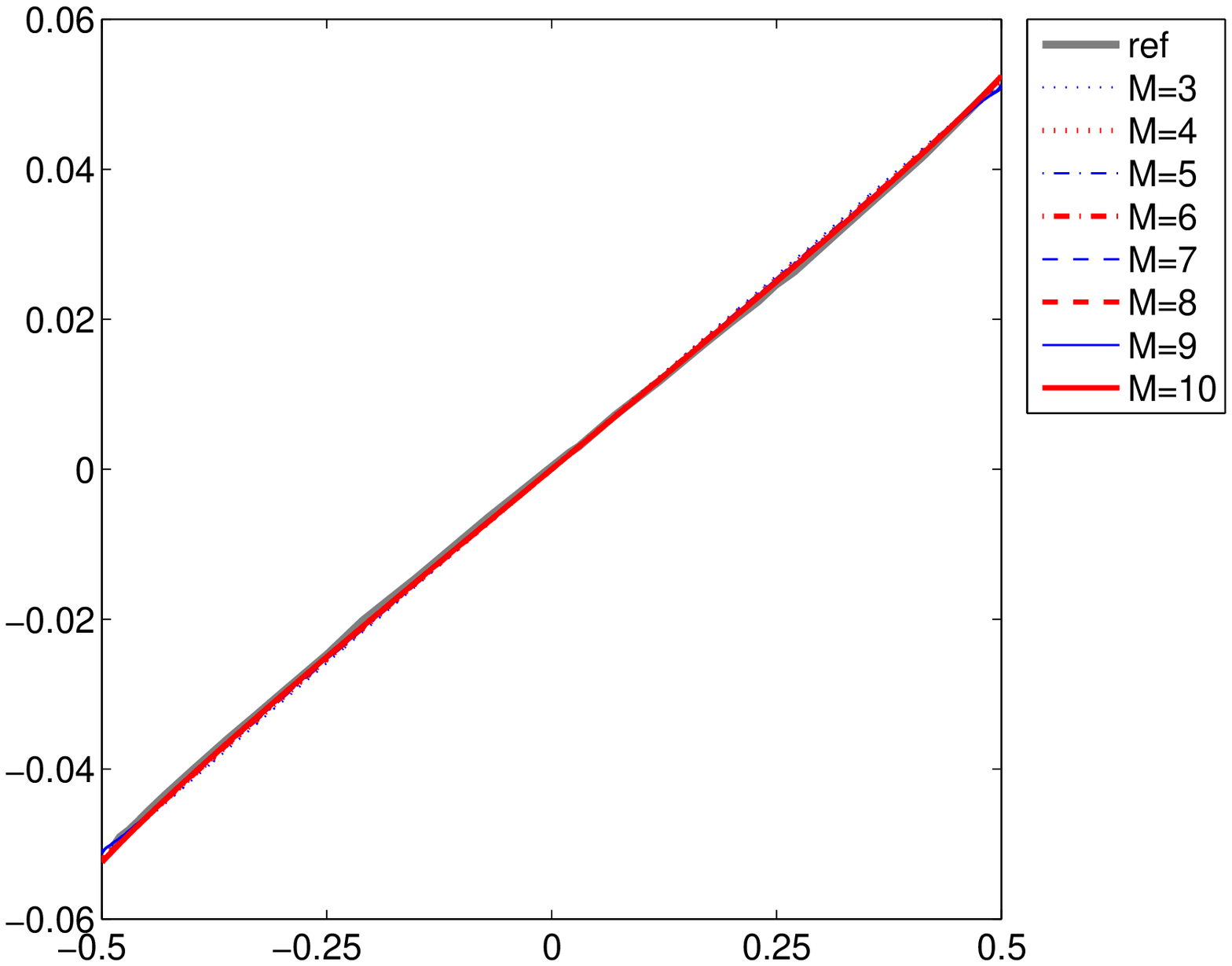}} }
\caption{Solution of the Couette flow for $\Kn=0.1199$ with $u^W =
  1.2577$ on a uniform grid of $N=2048$.}
  \label{fig:couette-Kn01-uw300}
\end{figure}

\begin{figure}[!htb]
{  \centering
  \subfigure[Density, $\rho$]{\includegraphics[width=0.231\textwidth]{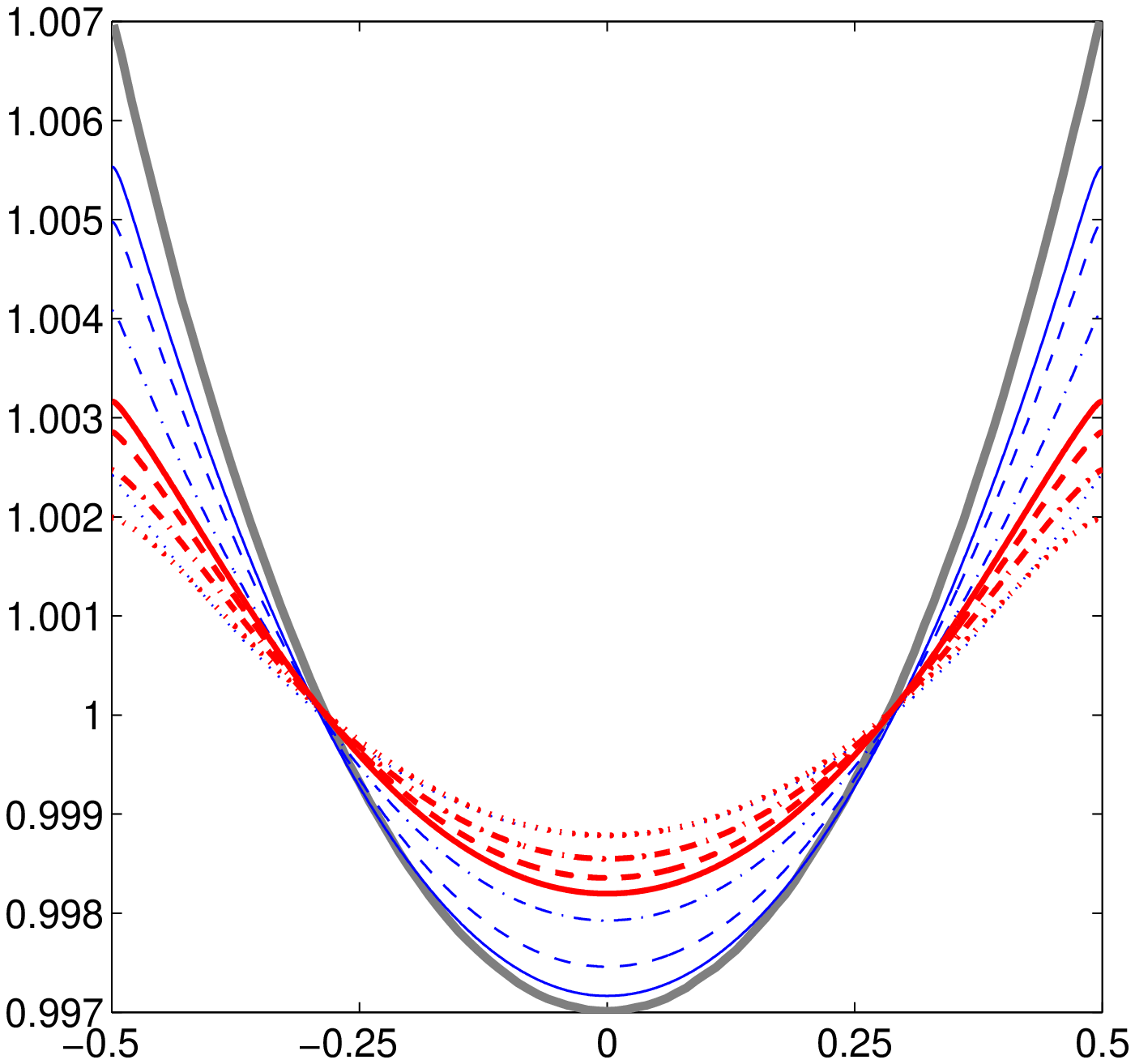}}
  \subfigure[Temperature, $\theta$]{\includegraphics[width=0.231\textwidth]{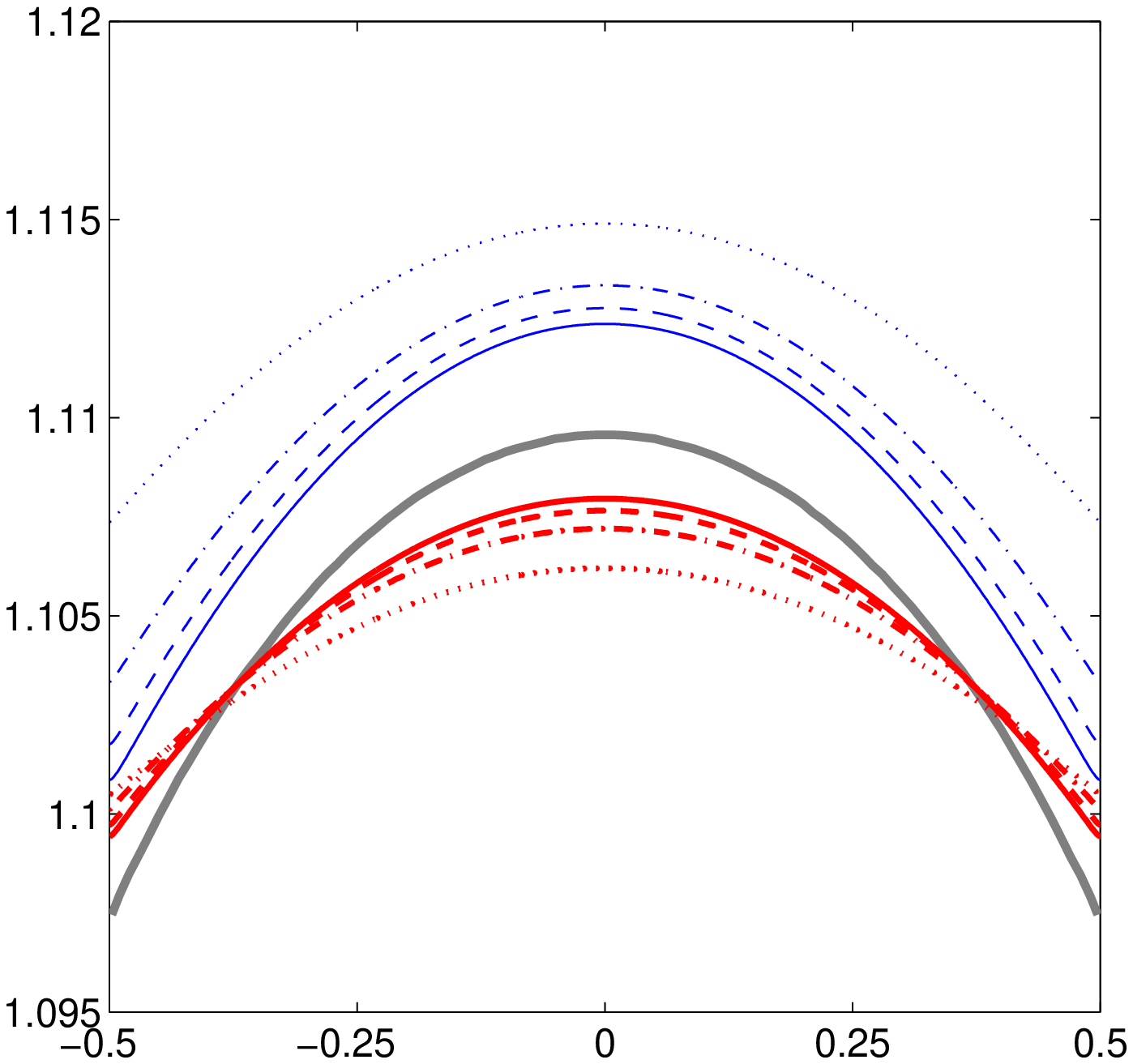}} 
  \subfigure[Shear stress, $\sigma_{12}$]{\includegraphics[width=0.236\textwidth]{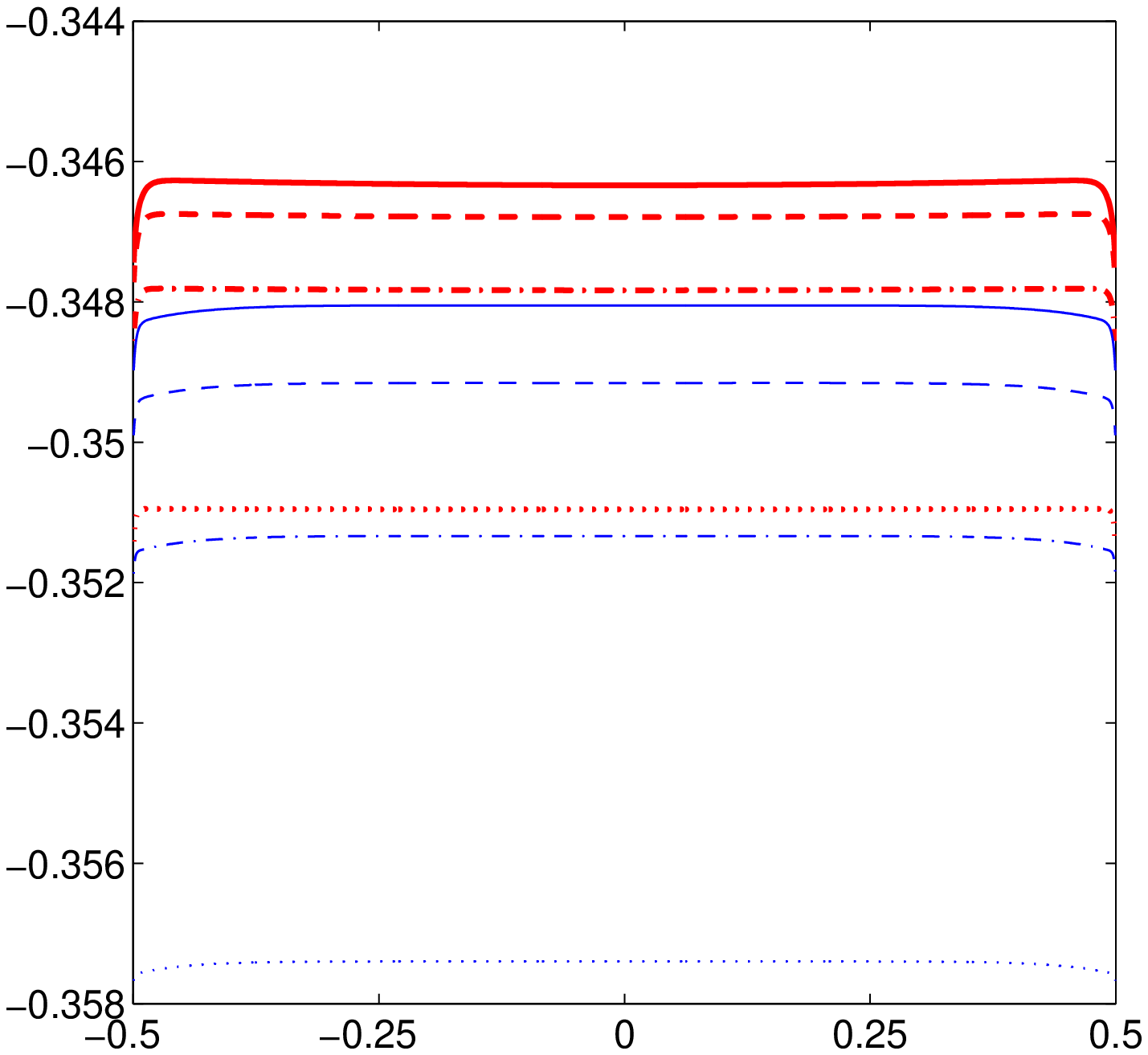}} 
  \subfigure[Heat flux, $q_1$]{\includegraphics[width=0.277\textwidth]{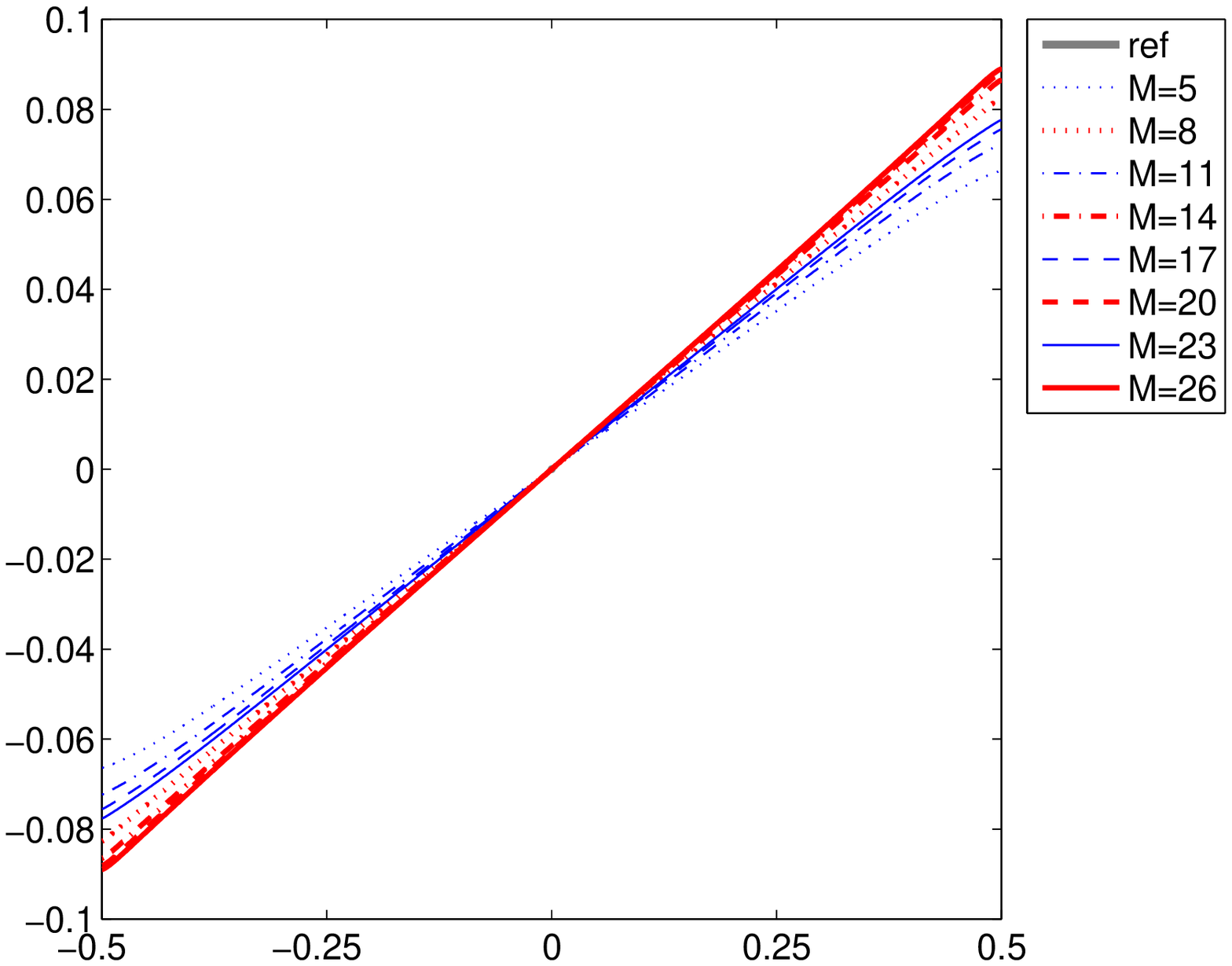}} }
\caption{Solution of the Couette flow for $\Kn=1.199$ with $u^W =
  1.2577$ on a uniform grid of $N=2048$.}
  \label{fig:couette-Kn10-uw300}
\end{figure}

We adopt the same settings as in \cite{Microflows1D, hu2014nmg}. To be
specific, the gas of argon is considered, and we have $\theta^W=1$,
$L=1$. The dimensionless collision frequency $\nu$ is given by
\begin{align}
  \label{eq:couette-nu}
  \nu = \sqrt{\frac{\pi}{2}} \frac{\Pr}{\Kn} \rho \theta^{1-w},
\end{align}
where $\Kn$ is the Knudsen number, and $w$ is the viscosity index. For
the gas of argon, the value of $w$ is $0.81$. With these parameters,
the proposed NMLM solver delivers exactly the steady-state solution
obtained in \cite{Microflows1D,hu2014nmg}. Since in
\cite{Microflows1D} the solution of the moment models has been
compared with the reference solution obtained in \cite{Mieussens2004},
and its convergence with respect to the order $M$ has been validated,
we omit any discussion on the accuracy and the convergence with
respect to $M$ of our solution. As examples, the steady-state solution
for $\Kn = 0.1199$ and $1.199$ with $u^W = 1.2577$ on a uniform grid
of $N = 2048$ are displayed in Figure \ref{fig:couette-Kn01-uw300} and
\ref{fig:couette-Kn10-uw300} respectively, in comparison to the
reference solution. It can be seen that the moment model of order
$M=10$ is enough to give satisfactory results for $\Kn=0.1199$, while
the moment model up to order $M=23$ or $26$ is still necessary to be
used for $\Kn=1.199$.


As pointed out in \cite{Microflows1D}, the moment models reduce
degrees of freedom significantly in comparison to the discrete
velocity method that was used in \cite{Mieussens2004}. While on the
other hand, we have observed from our computations that as a variation
of explicit time-integration scheme, the SGS-Richardson iteration
converges in general several times faster, consequently more
efficient, than the time-integration scheme employed in
\cite{Microflows1D}. Therefore, below we only investigate the
effectiveness of the multi-level strategy using lower-order model
correction to accelerate the convergence and the behavior of the
resulting NMLM solver. For comparison, all the computations start from
the same global equilibrium with
\begin{align}
  \label{eq:couette-initial}
  \rho^0 (x) = 1, \quad \bu^0(x) = 0, \quad \theta^0(x) = 1.
\end{align}

We perform the NMLM solver with different levels and order reduction
strategies for the moment model of various orders on three uniform
grids of $N=128$, $256$, $512$, respectively. Only some of numerical
results are shown in this paper, since the NMLM solver exhibits
similar features for all cases. In the tables given below, $K$ and $T$
represent respectively the total number and CPU seconds of the NMLM
iterations to reach the steady state, while $K_s$ and $T_s$ are
corresponding quantities of the single level solver.

First the Couette flow for $\Kn=0.1199$ and $u^W=1.2577$ is
considered. Table \ref{tab:couette-Kn01-uw300:nmm-M45} gives the
performance results for the case of order $M=4$ and $5$. The
corresponding convergence histories on the uniform grid of $N=512$ are
shown in Figure \ref{fig:couette-Kn01-uw300:nmm-res-history-N512-M4}
for $M=4$ and in Figure
\ref{fig:couette-Kn01-uw300:nmm-res-history-N512-M5} for $M=5$,
respectively. It is quite inspiring that the NMLM solver is effective
for such cases, where the order of the moment model is not very large.
For both cases, the convergence is accelerated and the total
computational cost, i.e., the CPU time, is reduced a lot, by the
multi-level NMLM iterations, in comparison to the single level
solver. It can be seen that for two-level NMLM solvers, the order
reduction strategy $m_{\sss l-1} = m_{\sss l} - 2$ converges faster
than the strategy $m_{\sss l-1} = m_{\sss l} - 1$. Moreover, the
computational cost of each NMLM iteration for the former strategy is
also less than the latter strategy, since the strategy $m_{\sss l-1} =
m_{\sss l} - 2$ employs a lower-order model correction with the order
less than the counterpart of the strategy $m_{\sss l-1} = m_{\sss l} -
1$. Thus, the overall performance of the strategy $m_{\sss l-1} =
m_{\sss l} - 2$ is better than the strategy $m_{\sss l-1} = m_{\sss l}
- 1$, when the same two levels is used in the NMLM solver. As the
total levels up to 3, the convergence rate of the NMLM solver becomes
better than both two-level NMLM solvers. Although the strategy
$m_{\sss l-1} = m_{\sss l} - 1$ becomes also more efficient as the
total levels increases, the three-level NMLM solver with the strategy
$m_{\sss l-1} = m_{\sss l} - 1$ would still not be more efficient than
the two-level NMLM solver with the strategy $m_{\sss l-1} = m_{\sss l}
- 2$. At last, it can also be found from Table
\ref{tab:couette-Kn01-uw300:nmm-M45} that the multi-level NMLM solver
behaves similar to the single level solver as well as the explicit
time-integration scheme. That is, the total number of NMLM iterations
doubles and the total CPU seconds quadruples, as the grid number $N$
doubles.

For the case of order $M=10$, the performance results are listed in
Table
\ref{tab:couette-Kn01-uw300:nmm-M10-p1}-\ref{tab:couette-Kn01-uw300:nmm-M10-p2},
and the corresponding convergence histories on the uniform grid of
$N=512$ are shown in Figure
\ref{fig:couette-Kn01-uw300:nmm-res-history-N512-M10}. Now the order
reduction strategy $m_{\sss l-1} = \lceil m_{\sss l} / 2 \rceil$ can
also be applied. It can be seen again that the multi-level NMLM
solvers for all three order reduction strategies could accelerate the
steady-state computation. In more details, when the NMLM solvers with
the same total levels are performed, the most efficient order
reduction strategy is $m_{\sss l-1} = \lceil m_{\sss l} / 2 \rceil$,
the second is $m_{\sss l-1} = m_{\sss l} - 2$, and the third is
$m_{\sss l-1} = m_{\sss l} - 1$, for that they are in descending sort
not only on the speed of convergence, but also on the computational
cost of each NMLM iteration. As the total levels increases, both the
convergence rate and the efficiency of the NMLM solver become better
for each order reduction strategy. However, the $8$-level NMLM solver
with the strategy $m_{\sss l-1} = m_{\sss l} - 1$ does still less
efficient than the $5$-level NMLM solver with the strategy $m_{\sss
  l-1} = m_{\sss l} - 2$, whereas the overall performance of the
latter solver is just close to the $3$-level NMLM solver with the
strategy $m_{\sss l-1} = \lceil m_{\sss l} / 2 \rceil$, for which the
total computational cost is saved by approximately more than $80\%$ in
comparison to the single level solver. In addition, we have again that
the total number of NMLM iterations doubles and the total CPU seconds
quadruples, as the grid number $N$ doubles.

As mentioned previous, the moment model up to order $M=23$ or $26$
should be taken into consideration when $\Kn=1.199$. A partial
performance results are shown in Table
\ref{tab:couette-Kn10-uw300:nmm-M23-p2} for the case of order $M=23$,
and in Table \ref{tab:couette-Kn10-uw300:nmm-M26-p2} for the case of
order $M=26$, respectively. The corresponding convergence histories
are plotted in Figure
\ref{fig:couette-Kn10-uw300:nmm-res-history-N512-M23}-\ref{fig:couette-Kn10-uw300:nmm-res-history-N512-M26}.
We do not present results of the NMLM solver with the strategy
$m_{\sss l-1} = m_{\sss l} - 1$ here, since compared with the single
level solver it turns out a little improvement in efficiency, although
the speed of convergence is raised much. This is reasonable by noting
that the order of the lower-order problem just reduces $1$ at each
level, for example, the order sequence for $M=26$ is $26, 25, 24,
\ldots$, which indicates that the lower-order model correction still
takes a lot of computational cost.  In fact, the computational cost of
lower-order model correction can not be underestimated, even when the
strategy $m_{\sss l-1} = m_{\sss l} - 2$, giving the order sequence
$26,24,22,\ldots$ for $M=26$, is adopted. Moreover, it can be seen
that the multi-level NMLM solver has some degeneracy, especially for
the solver with the strategy $m_{\sss l-1} = \lceil m_{\sss l} / 2
\rceil$. As a result, the overall performance of the multi-level NMLM
solver would not be as good as the solver when $\Kn=0.1199$, although
the efficiency is still improved much compared with the single level
solver. Furthermore, unlike the observation when $\Kn=0.1199$, the
convergence rate of the strategy $m_{\sss l-1} = \lceil m_{\sss l} / 2
\rceil$ is worse than the strategy $m_{\sss l-1} = m_{\sss l} -
2$. However, with the help of great reduction of the computational
cost at each NMLM iteration, the strategy $m_{\sss l-1} = \lceil
m_{\sss l} / 2 \rceil$ finally exhibits more efficient than the
strategy $m_{\sss l-1} = m_{\sss l} - 2$. On the other hand,
oscillations of the residual are now observed at the beginning
iterations of single level solver. For the multi-level NMLM solvers,
the oscillations become more severer, and may introduce
instability of the solver. Actually, the $5$-level NMLM solver with
the strategy $m_{\sss l-1} = \lceil m_{\sss l} / 2 \rceil$ breaks down
in our computations. In view of these, a possible way of taking both
efficiency and stability into account might be to adopt the order
reduction strategy $m_{\sss l-1} = m_{\sss l} - \delta m$, such that
$\delta m > 2$ and $m_{\sss l} - \delta m > \lceil m_{\sss l} / 2
\rceil $. At last, we have again that the convergence rate is improved
by the multi-level NMLM solver as the total levels increases, and the
multi-level NMLM solvers behave similarly to the single level solver,
as the grid number $N$ doubles.

\begin{figure}[!htb]
  \centering
  \subfigure[$\Kn=0.1199$]{\includegraphics[width=0.49\textwidth]{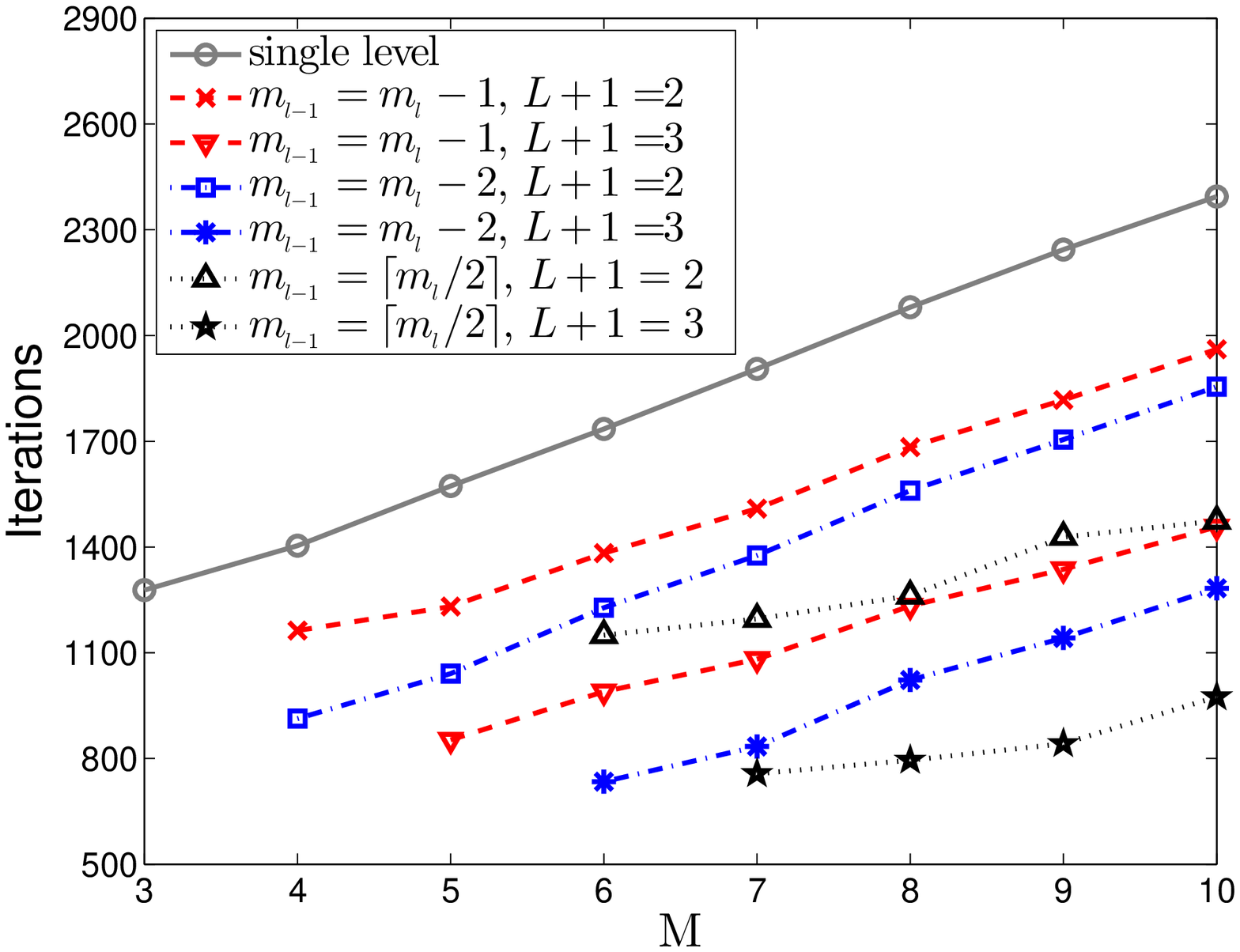}}
  \subfigure[$\Kn=1.199$]{\includegraphics[width=0.49\textwidth]{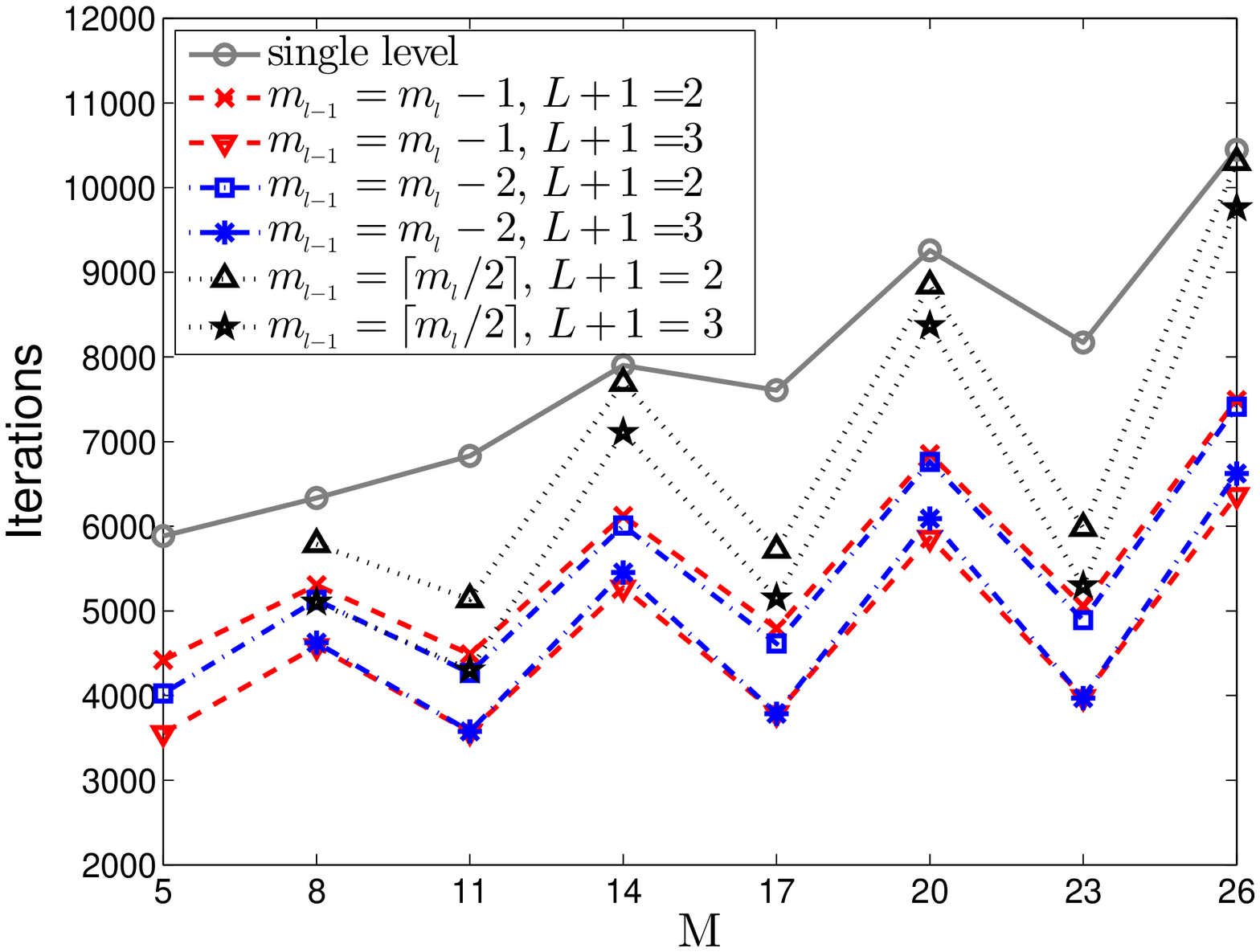}}
  \caption{Total iterations in terms of $M$ of the NMLM solver for the
    Couette flow with $u^W = 1.2577$ on a uniform grid of $N=512$. The
    total iterations of the single level solver is rescaled by a
    factor of $12$ for $\Kn=0.1199$ and $8$ for $\Kn=1.199$
    respectively.}
  \label{fig:couette:nmm-iters-orders}
\end{figure}

It is noted from Table
  \ref{tab:couette-Kn10-uw300:nmm-M23-p2}-\ref{tab:couette-Kn10-uw300:nmm-M26-p2}
  that the total iterations $K$ is almost doubled as $M$ increases
  from $23$ to $26$, while the total iterations $K$ shown in Table
  \ref{tab:couette-Kn01-uw300:nmm-M45}-\ref{tab:couette-Kn01-uw300:nmm-M10-p2}
  increases much slower as $M$ increases from $4$ to $10$. The
  significant difference is mainly due to the different performance of
  the smoothing operator (equivalently the single level solver) with
  respect to the Knudsen number $\Kn$. To see it in more detail, we
  plot $K$ in terms of $M$ for the NMLM solver in Figure
  \ref{fig:couette:nmm-iters-orders}. It can be seen that the total
  iterations $K$ of the single level solver increases linearly with
  respect to $M$ for the case $\Kn=0.1199$, whereas for the case
  $\Kn=1.199$ the total iterations $K$ of the single level solver
  shows a strong difference with respect to the parity of $M$,
  especially for a larger $M$. To be specific, the total iterations
  $K$ increases linearly with a smaller rate with respect to odd $M$,
  and with a larger rate with respect to even $M$. As for the
  two-level and three-level NMLM solvers with the same order reduction
  strategy, we can see the total iterations $K$ increases linearly
  with similar rate with respect to $M$, in comparison to the
  corresponding single level solver. 

In summary, it is effective to accelerate the steady-state computation
by using the multi-level NMLM solver. The convergence rate would
become better as the total levels increases, and the total
computational cost is then saved a lot by comparing with the single
level solver. Among three order reduction strategies, the strategy
$m_{\sss l-1} = \lceil m_{\sss l} / 2 \rceil$ would be most efficient,
followed with the strategy $m_{\sss l-1} = m_{\sss l} - 2$, and then
the strategy $m_{\sss l-1} = m_{\sss l} - 1$.

\begin{table}[!ht]
  \centering\scriptsize
  \begin{tabular}{c|c||c|c|c||c|cc|c}
    \hline\hline
    \multicolumn{2}{c||}{} & \multicolumn{3}{c||}{$M=4$} & \multicolumn{4}{c}{$M=5$}\\
    \hline
    \multicolumn{2}{c||}{} & & \multicolumn{1}{c|}{$m_{\sss l-1} = m_{\sss l} - 1$} & \multicolumn{1}{c||}{$m_{\sss l-1} = m_{\sss l} - 2$} & &\multicolumn{2}{c|}{$m_{\sss l-1} = m_{\sss l} - 1$} & \multicolumn{1}{c}{$m_{\sss l-1} = m_{\sss l} - 2$}\\
    \hline
    \multicolumn{2}{c||}{$L+1$} & 1 & 2 & 2 & 1 & 2 & 3 & 2 \\
    \hline\hline
    \multirow{4}{*}{\begin{sideways}$N=128$\end{sideways}} 
    & $K$          &     3955&     263&         210&       4217&      261&         174&         228 \\
    & $T$       &  155.558&  80.320&      49.219&    243.090&  122.534&      59.466&      84.243 \\
    & $K_s/K$      &      1.0&  15.038&      18.833&        1.0&   16.157&      24.236&      18.496 \\
    & $T_s/T$      &      1.0&   1.937&       3.161&        1.0&    1.984&       4.088&       2.886 \\
    \hline\hline
    \multirow{4}{*}{\begin{sideways}$N=256$\end{sideways}} 
    & $K$          &     8178&      557&         440&      9064&     591&         410&         500 \\
    & $T$       &  633.756&  348.631&     209.665&  1059.542& 573.913&     420.866&     371.722 \\
    & $K_s/K$      &      1.0&   14.682&      18.586&       1.0&  15.337&      22.107&      18.128 \\
    & $T_s/T$      &      1.0&    1.818&       3.023&       1.0&   1.846&       2.518&       2.850 \\
    \hline\hline
    \multirow{4}{*}{\begin{sideways}$N=512$\end{sideways}} 
    & $K$          &    16848&     1163&         913&    18875&     1231&         853&        1041 \\
    & $T$       & 2390.169& 1054.177&     871.586& 4433.716& 2350.401&    1677.407&    1501.216 \\
    & $K_s/K$      &      1.0&   14.487&      18.453&      1.0&   15.333&      22.128&      18.132 \\
    & $T_s/T$      &      1.0&    2.267&       2.742&      1.0&    1.886&       2.643&       2.953 \\
    \hline\hline
  \end{tabular}
  \caption{Performance of the NMLM solver for the Couette flow with $\Kn=0.1199$, $u^W = 1.2577$ and $M=4,5$.}
  \label{tab:couette-Kn01-uw300:nmm-M45}
\end{table}

\begin{table}[!ht]
  \centering\footnotesize
  \begin{tabular}{c|c||ccccccc}
    \hline\hline
    \multicolumn{2}{c||}{} & \multicolumn{6}{c}{$m_{\sss l-1} = m_{\sss l} - 1$}\\
    \hline
    \multicolumn{2}{c||}{$L+1$} & 2 & 3 & 4 & 5 & 6 & 7 & 8 \\
    \hline\hline
    \multirow{4}{*}{\begin{sideways}$N=128$\end{sideways}} 
    & $K$          &       476&         354&         276&         224&         189&         168&         156\\
    & $T$       &  1287.045&    1012.705&     908.478&     773.019&     681.257&     617.487&     573.095\\
    & $K_s/K$      &    14.639&      19.684&      25.246&      31.107&      36.868&      41.476&      44.667\\
    & $T_s/T$      &     1.481&       1.882&       2.097&       2.465&       2.797&       3.086&       3.325\\
    \hline\hline
    \multirow{4}{*}{\begin{sideways}$N=256$\end{sideways}} 
    & $K$          &       963&         716&         560&         454&         380&         332&         303\\
    & $T$       &  4693.688&    3787.012&    3767.613&    3176.846&    2716.522&    2189.956&    2216.659\\
    & $K_s/K$      &    14.652&      19.707&      25.196&      31.079&      37.132&      42.500&      46.568\\
    & $T_s/T$      &     1.810&       2.243&       2.255&       2.674&       3.127&       3.879&       3.833\\
    \hline\hline
    \multirow{4}{*}{\begin{sideways}$N=512$\end{sideways}} 
    & $K$          &      1961&        1457&        1139&         921&         768&         663&         597\\
    & $T$       & 16760.802&   16049.136&   10434.860&   12920.731&   10710.200&    9641.965&    6954.572\\
    & $K_s/K$      &    14.651&      19.719&      25.224&      31.194&      37.409&      43.333&      48.124\\
    & $T_s/T$      &     1.681&       1.756&       2.700&       2.181&       2.631&       2.922&       4.051\\
    \hline\hline
  \end{tabular}
  \caption{Performance of the NMLM solver for the Couette flow with $\Kn=0.1199$, $u^W = 1.2577$ and $M=10$ (part I).}
  \label{tab:couette-Kn01-uw300:nmm-M10-p1}
\end{table}

\begin{table}[!ht]
  \centering\footnotesize
  \begin{tabular}{c|c||cccc|cc|c}
    \hline\hline
    \multicolumn{2}{c||}{} & \multicolumn{4}{c|}{$m_{\sss l-1} = m_{\sss l} - 2$} & \multicolumn{2}{c|}{$m_{\sss l-1} = \lceil m_{\sss l} / 2 \rceil$} & \\
    \hline
    \multicolumn{2}{c||}{$L+1$} & 2 & 3 & 4 & 5 & 2 & 3 & 1\\
    \hline\hline
    \multirow{4}{*}{\begin{sideways}$N=128$\end{sideways}} 
    & $K$          &        450&         312&         221&         165&         357&         237&      6968\\
    & $T$       &    789.648&     597.272&     497.074&     363.046&     531.005&     334.993&  1905.477\\
    & $K_s/K$      &     15.484&      22.333&      31.529&      42.230&      19.518&      29.401&       1.0\\
    & $T_s/T$      &      2.413&       3.190&       3.833&       5.249&       3.588&       5.688&       1.0\\
    \hline\hline
    \multirow{4}{*}{\begin{sideways}$N=256$\end{sideways}} 
    & $K$          &        912&         631&         446&         327&         724&         479&     14110\\
    & $T$       &   4233.031&    2572.533&    2004.516&    1445.255&    1951.228&    1356.081&  8495.687\\
    & $K_s/K$      &     15.471&      22.361&      31.637&      43.150&      19.489&      29.457&       1.0\\
    & $T_s/T$      &      2.007&       3.302&       4.238&       5.878&       4.354&       6.265&       1.0\\
    \hline\hline
    \multirow{4}{*}{\begin{sideways}$N=512$\end{sideways}} 
    & $K$          &       1855&        1283&         903&         651&        1474&         975&     28730\\
    & $T$       &  13206.382&   11231.685&    8146.139&    4405.231&    7564.141&    5498.801& 28174.869\\
    & $K_s/K$      &     15.488&      22.393&      31.816&      44.132&      19.491&      29.467&       1.0\\
    & $T_s/T$      &      2.133&       2.509&       3.459&       6.396&       3.725&       5.124&       1.0\\
    \hline\hline
  \end{tabular}
  \caption{Performance of the NMLM solver for the Couette flow with $\Kn=0.1199$, $u^W = 1.2577$ and $M=10$ (part II).}
  \label{tab:couette-Kn01-uw300:nmm-M10-p2}
\end{table}

\begin{table}[!ht]
\hspace*{-2.5em}
  \centering\scriptsize
  \begin{tabular}{c|c||ccccc|ccc}
    \hline\hline
    \multicolumn{2}{c||}{} & \multicolumn{5}{c|}{$m_{\sss l-1} = m_{\sss l} - 2$} & \multicolumn{3}{c}{$m_{\sss l-1} = \lceil m_{\sss l} / 2 \rceil$} \\
    \hline
    \multicolumn{2}{c||}{$L+1$} & 4 & 5 & 6 & 7 & 8 & 2 & 3 & 4 \\
    \hline\hline
    \multirow{4}{*}{\begin{sideways}$N=128$\end{sideways}} 
    & $K$          &         796&         660&         481&         506&         530&        1440&        1440&        1338\\
    & $T$       &   26508.912&   22966.565&   17201.641&   17694.056&   13069.990&   18894.415&   19519.867&   10330.927\\
    & $K_s/K$      &      16.991&      20.492&      28.119&      26.729&      25.519&       9.392&       9.392&      10.108\\
    & $T_s/T$      &       1.165&       1.345&       1.796&       1.746&       2.364&       1.635&       1.583&       2.991\\
    \hline\hline
    \multirow{4}{*}{\begin{sideways}$N=256$\end{sideways}} 
    & $K$          &        1593&        1386&        1220&        1099&        1007&        2405&        2602&        2409\\
    & $T$       &  105550.167&   96564.039&   86944.937&   79480.186&   57305.984&   62363.809&   72069.420&   65186.066\\
    & $K_s/K$      &      16.559&      19.032&      21.621&      24.002&      26.195&      10.968&      10.138&      10.950\\
    & $T_s/T$      &       1.396&       1.526&       1.694&       1.854&       2.571&       2.362&       2.044&       2.260\\
    \hline\hline
    \multirow{4}{*}{\begin{sideways}$N=512$\end{sideways}} 
    & $K$          &        3392&        2998&        2717&        2511&        2353&        5979&        5295&        5074\\
    & $T$       &  408209.027&  363005.320&  366497.127&  316449.369&  268045.359&  300242.049&  266417.091&  216497.893\\
    & $K_s/K$      &      19.268&      21.801&      24.055&      26.029&      27.776&      10.931&      12.343&      12.881\\
    & $T_s/T$      &       1.903&       2.140&       2.119&       2.455&       2.898&       2.587&       2.915&       3.588\\
    \hline\hline
  \end{tabular}
  \caption{Performance of the NMLM solver for the Couette flow with $\Kn=1.199$, $u^W = 1.2577$ and $M=23$.}
  \label{tab:couette-Kn10-uw300:nmm-M23-p2}
\end{table}

\begin{table}[!ht]
\hspace*{-2.5em}
  \centering\scriptsize
  \begin{tabular}{c|c||ccccc|ccc}
    \hline\hline
    \multicolumn{2}{c||}{} & \multicolumn{5}{c|}{$m_{\sss l-1} = m_{\sss l} - 2$} & \multicolumn{3}{c}{$m_{\sss l-1} = \lceil m_{\sss l} / 2 \rceil$} \\
    \hline
    \multicolumn{2}{c||}{$L+1$} & 4 & 5 & 6 & 7 & 8 & 2 & 3 & 4 \\
    \hline\hline
    \multirow{4}{*}{\begin{sideways}$N=128$\end{sideways}} 
    & $K$          &        1559&        1472&        1410&        1363&        1326&        2627&        2496&        2405\\ 
    & $T$       &   67950.969&   73383.501&   58109.432&   76334.248&   68269.101&   48025.102&   47618.776&   46879.272\\ 
    & $K_s/K$      &      13.528&      14.327&      14.957&      15.473&      15.905&       8.028&       8.450&       8.769\\ 
    & $T_s/T$      &       1.124&       1.041&       1.315&       1.001&       1.119&       1.591&       1.605&       1.630\\ 
    \hline\hline
    \multirow{4}{*}{\begin{sideways}$N=256$\end{sideways}} 
    & $K$          &        3083&        2911&        2789&        2696&        2622&        5190&        4920&        4715\\ 
    & $T$       &  276985.452&  297746.135&  296101.816&  275546.571&  212048.261&  211041.390&  194471.039&  178270.975\\
    & $K_s/K$      &      13.586&      14.389&      15.019&      15.537&      15.975&       8.071&       8.514&       8.884\\ 
    & $T_s/T$      &       1.265&       1.177&       1.183&       1.271&       1.652&       1.660&       1.801&       1.965\\ 
    \hline\hline
    \multirow{4}{*}{\begin{sideways}$N=512$\end{sideways}} 
    & $K$          &        6116&        5778&        5536&        5354&       5207&       10303&        9761&        9319\\ 
    & $T$       & 1147507.284&  887609.149&  992797.036&  889783.318& 801075.758&  607645.946&  663480.962&  555586.099\\ 
    & $K_s/K$      &      13.664&      14.463&      15.095&      15.609&     16.049&       8.111&       8.561&       8.967\\ 
    & $T_s/T$      &       1.158&       1.497&       1.338&       1.493&      1.658&       2.186&       2.002&       2.391\\ 
    \hline\hline
  \end{tabular}
  \caption{Performance of the NMLM solver for the Couette flow with $\Kn=1.199$, $u^W = 1.2577$ and $M=26$.}
  \label{tab:couette-Kn10-uw300:nmm-M26-p2}
\end{table}

\newcommand\drawCouetteNMMHistory[6]{\begin{figure}[!htb]
  \centering
  {\includegraphics[width=0.5\textwidth]{couette_#1_#3_nmm_res_iters_N#5_M#6.eps}}\hfill
  {\includegraphics[width=0.5\textwidth]{couette_#1_#3_nmm_res_cputime_N#5_M#6.eps}}
  \caption{Convergence history of the NMLM solver for the Couette flow
    with $\Kn=#2$, $u^W = #4$ and $M=#6$ on a uniform grid of $N=#5$.}
  \label{fig:couette-#1-#3:nmm-res-history-N#5-M#6}
\end{figure}
}
\drawCouetteNMMHistory{Kn01}{0.1199}{uw300}{1.2577}{512}{4}
\drawCouetteNMMHistory{Kn01}{0.1199}{uw300}{1.2577}{512}{5}
\drawCouetteNMMHistory{Kn01}{0.1199}{uw300}{1.2577}{512}{10}
\drawCouetteNMMHistory{Kn10}{1.199}{uw300}{1.2577}{512}{23}
\drawCouetteNMMHistory{Kn10}{1.199}{uw300}{1.2577}{512}{26}

\subsection{The force driven Poiseuille flow}
\label{sec:num-ex-poiseuille}
\begin{figure}[!htb]
{  \centering
  \subfigure[Density, $\rho$]{\includegraphics[width=0.233\textwidth]{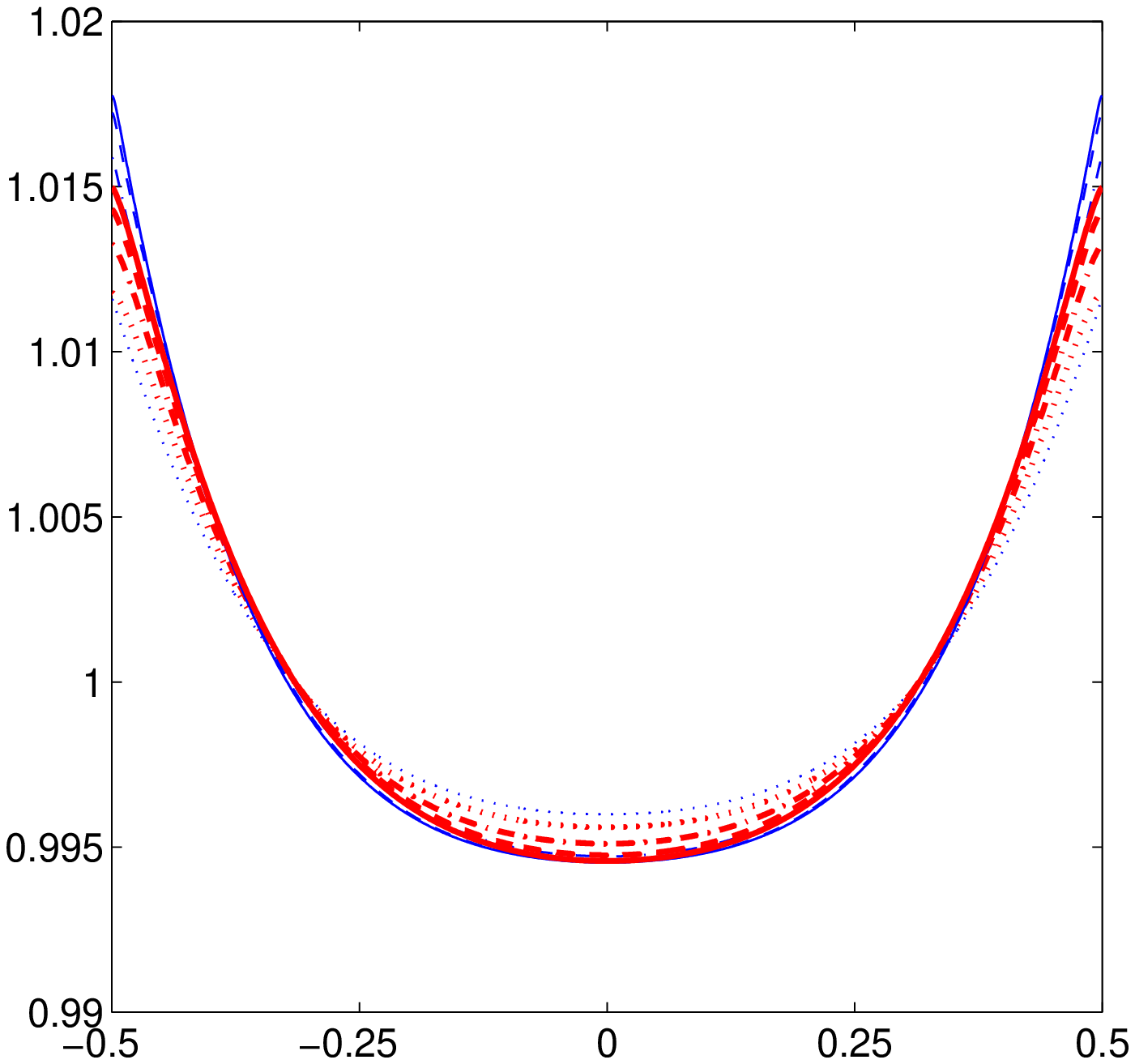}}
  \subfigure[Temperature, $\theta$]{\includegraphics[width=0.233\textwidth]{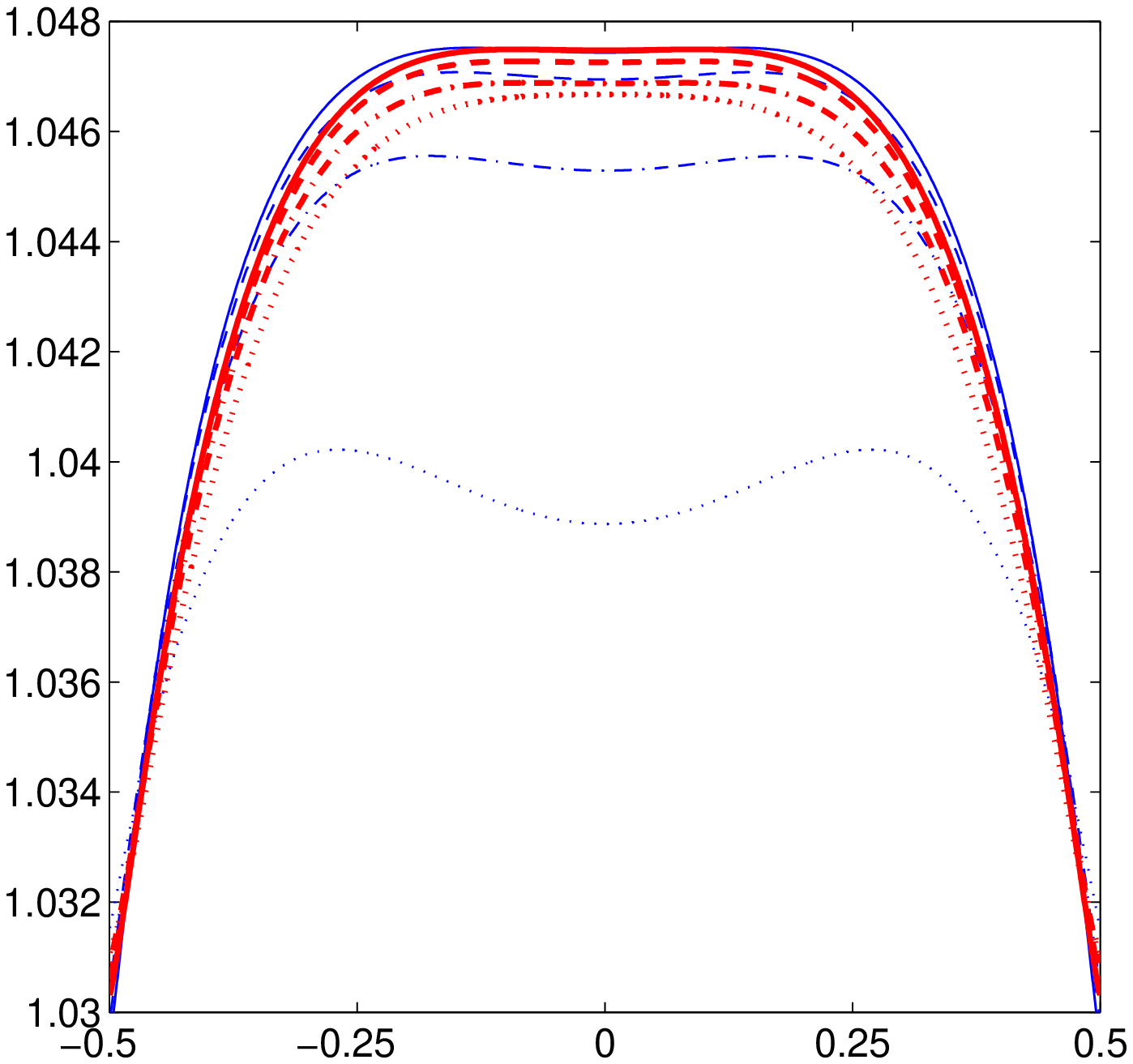}} 
  \subfigure[Normal stress, $\sigma_{11}$]{\includegraphics[width=0.228\textwidth]{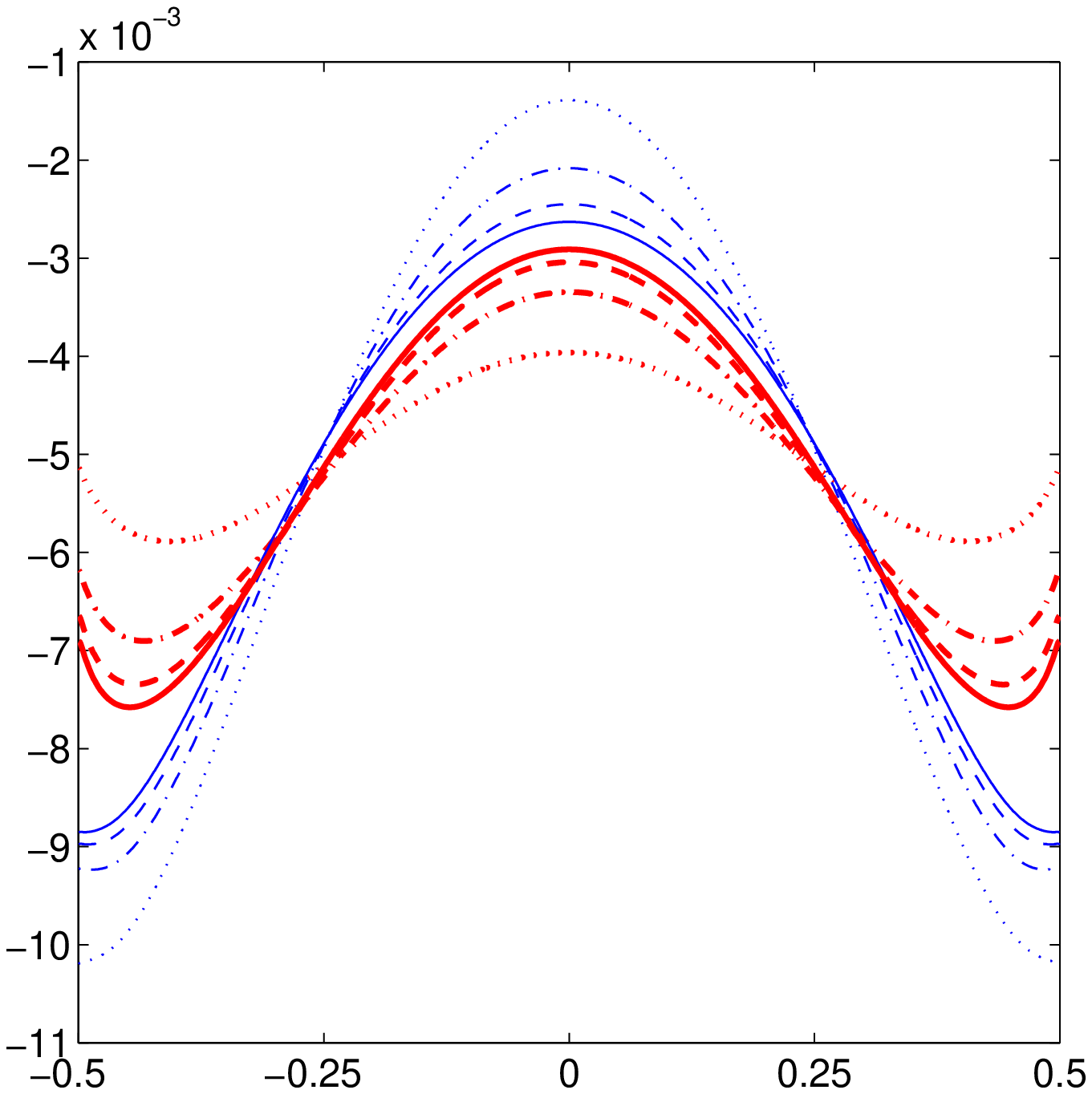}} 
  \subfigure[Heat flux, $q_2$]{\includegraphics[width=0.279\textwidth]{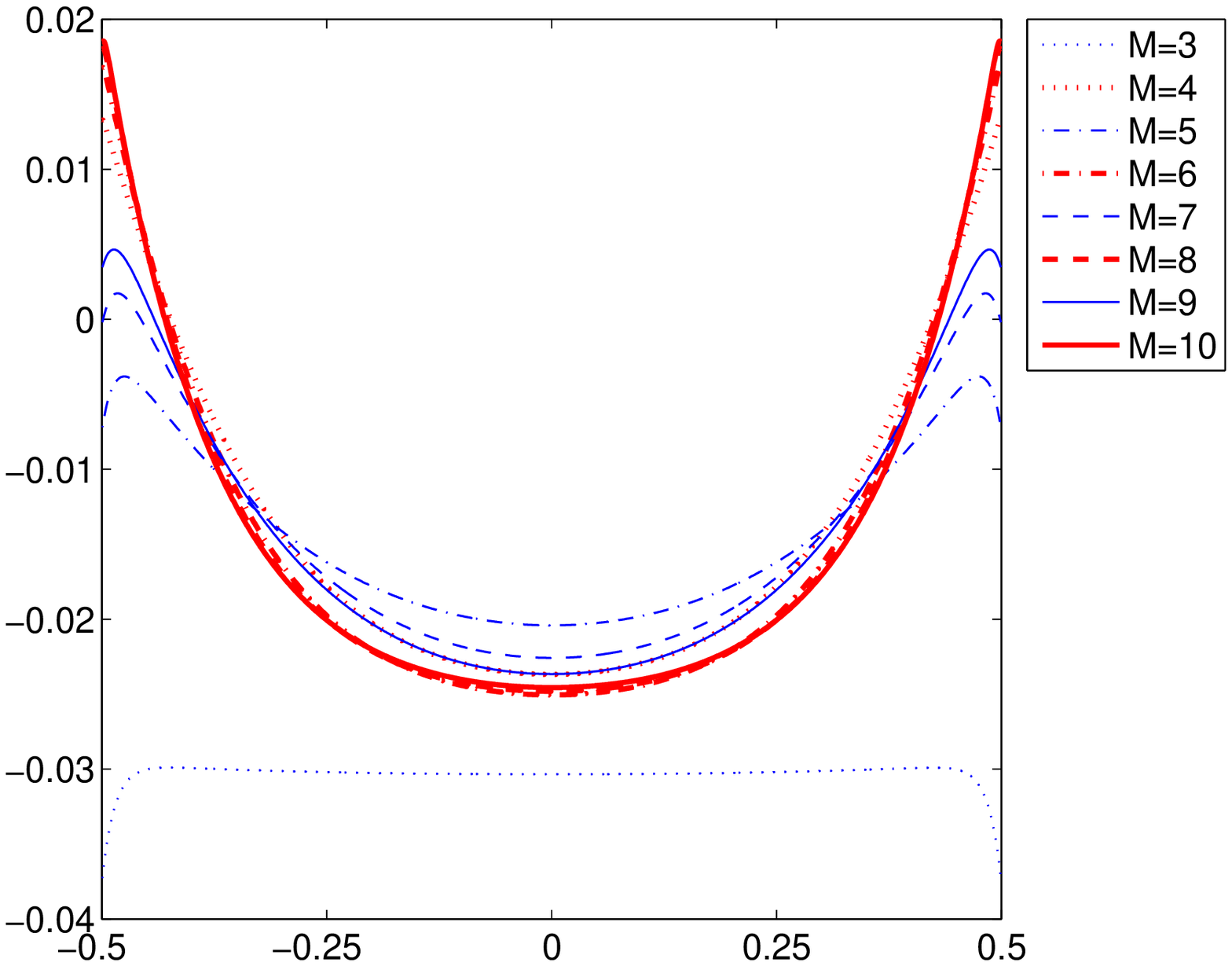}} }
\caption{Solution of the force driven Poiseuille flow on a uniform
  grid of $N=2048$.}
  \label{fig:poiseuille-sol}
\end{figure}

The force driven Poiseuille flow is another benchmark test frequently
investigated in the literatures \cite{Garcia, Li, Xu2007,
  hu2014nmg}. Similar to the Couette flow, there are two infinite
parallel plates, which are separated by a distance of $L=1$, and have
the same temperature of $\theta^W=1$. However, both plates are
stationary now, and the gas between them is driven by an external
constant force, which is set as $\bF=(0, 0.2555,0)^T$ in our
tests. Additionally, the collision frequency $\nu$ is given by the
hard sphere model as
\begin{align}
  \label{eq:poiseuille-nu}
  \nu = \frac{16}{5}\sqrt{\frac{\theta}{2 \pi}} \frac{\Pr}{\Kn} \rho,
\end{align}
and the Knudsen number $\Kn = 0.1$ is considered. With these settings,
the steady-state solution obtained by the NMLM solver is shown in
Figure \ref{fig:poiseuille-sol}, which recovers exactly the
steady-state solution presented in \cite{hu2014nmg}.

\begin{figure}[!htb]
  \centering
  {\includegraphics[width=0.5\textwidth]{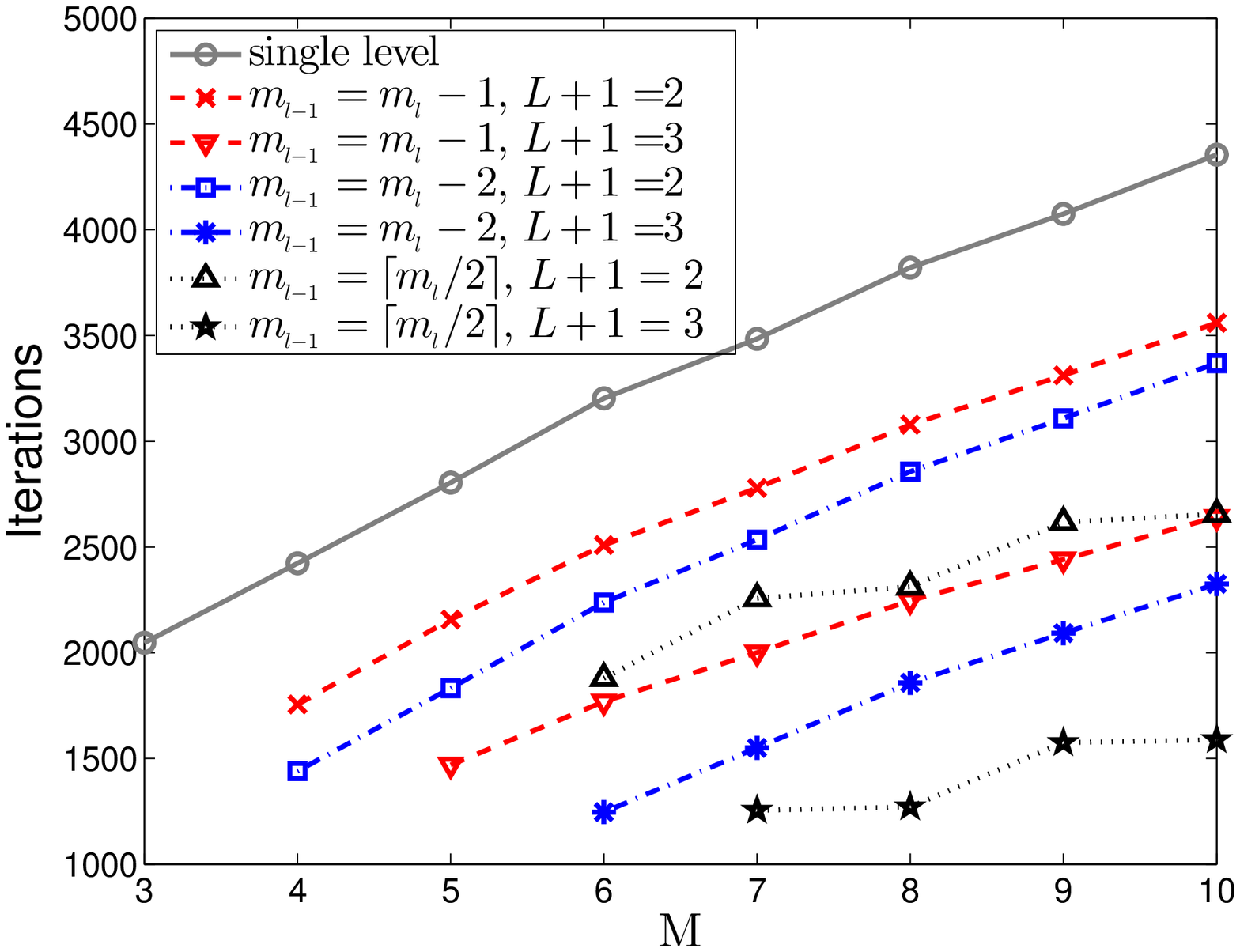}}
  \caption{Total iterations in terms of $M$ of the NMLM solver for the
    Poiseuille flow on a uniform grid of $N=512$. The total iterations
    of the single level solver is rescaled by a factor of $12$.}
  \label{fig:poiseuille:nmm-iters-orders}
\end{figure}

We still omit the discussion on the accuracy and the convergence of
the solution with respect to order $M$, and focus on the behavior of
the proposed NMLM solver. As the Couette flow, the NMLM solvers, with
different levels and order reduction strategies for the moment model
of various orders on three uniform grids of $N=128$, $256$, $512$, are
performed. The computations also begin with the global equilibrium
\eqref{eq:couette-initial}. Again just partial numerical results are
shown here, for similar features can be observed for all cases. To be
specific, the performance results are given in Table
\ref{tab:poiseuille:nmm-M45} for the case of order $M=4$, $5$, and in
Table \ref{tab:poiseuille:nmm-M10-p1}-\ref{tab:poiseuille:nmm-M10-p2}
for the case of order $M=10$, respectively. The corresponding
convergence histories on the uniform grid of $N=512$ are displayed
respectively in Figure \ref{fig:poiseuille:nmm-res-history-N512-M4}
for $M=4$, in Figure \ref{fig:poiseuille:nmm-res-history-N512-M5} for
$M=5$, and in Figure \ref{fig:poiseuille:nmm-res-history-N512-M10} for
$M=10$. The total iterations $K$ in terms of $M$ for the NMLM
  solver is presented in Figure
  \ref{fig:poiseuille:nmm-iters-orders}. All these results show that
the multi-level NMLM solver is able to accelerate the steady-state
computation significantly.

In comparison to results of the Couette flow with $\Kn=0.1199$, a
similar behavior of the multi-level NMLM solver can be observed. In
more details, we can see that the most efficient order reduction
strategy is $m_{\sss l-1} = \lceil m_{\sss l} / 2 \rceil$, the second
is $m_{\sss l-1} = m_{\sss l} - 2$, and the third is $m_{\sss l-1} =
m_{\sss l} - 1$. As can be seen from the tables, the ratio of $K_s$
and $K$ are all consistent with those for the Couette
flow. Consequently, the convergence rate of the multi-level NMLM
solver with all three order reduction strategies increase as the total
levels increases, and the total computational cost is saved greatly in
comparison to the single level solver. In addition, as the grid number
$N$ doubles, all multi-level NMLM solver show similar features as the
single level solver. Thus, the acceleration ratio will be maintained
even when a more fine spatial grid is adopted.

\begin{table}[!ht]
  \centering\scriptsize
  \begin{tabular}{c|c||c|c|c||c|cc|c}
    \hline\hline
    \multicolumn{2}{c||}{} & \multicolumn{3}{c||}{$M=4$} & \multicolumn{4}{c}{$M=5$}\\
    \hline
    \multicolumn{2}{c||}{} & & \multicolumn{1}{c|}{$m_{\sss l-1} = m_{\sss l} - 1$} & \multicolumn{1}{c||}{$m_{\sss l-1} = m_{\sss l} - 2$} & & \multicolumn{2}{c|}{$m_{\sss l-1} = m_{\sss l} - 1$} & \multicolumn{1}{c}{$m_{\sss l-1} = m_{\sss l} - 2$}\\ 
    \hline
    \multicolumn{2}{c||}{$L+1$} & 1 & 2 & 2 & 1 & 2 & 3 & 2 \\
    \hline\hline
    \multirow{4}{*}{\begin{sideways}$N=128$\end{sideways}} 
    & $K$          &      6660&        405&         330&      7627&        490&         334&         416\\
    & $T$       &   168.395&     73.819&      80.196&   300.558&    140.285&     156.427&     111.988\\
    & $K_s/K$      &       1.0&     16.444&      20.182&       1.0&     15.565&      22.835&      18.334\\
    & $T_s/T$      &       1.0&      2.281&       2.100&       1.0&      2.142&       1.921&       2.684\\
    \hline\hline
    \multirow{4}{*}{\begin{sideways}$N=256$\end{sideways}} 
    & $K$          &     14111&        855&         699&     16219&       1040&         709&         883\\
    & $T$       &   729.814&    502.633&     338.805&  1270.682&    837.806&     757.453&     677.226\\
    & $K_s/K$      &       1.0&     16.504&      20.187&       1.0&     15.595&      22.876&      18.368\\
    & $T_s/T$      &       1.0&      1.452&       2.154&       1.0&      1.517&       1.678&       1.876\\
    \hline\hline
    \multirow{4}{*}{\begin{sideways}$N=512$\end{sideways}} 
    & $K$          &     29077&       1756&        1441&     33653&       2157&        1470&        1832\\
    & $T$       &  2915.750&   2113.058&    1094.692&  6382.696&   3276.133&    3036.576&    2181.178\\
    & $K_s/K$      &       1.0&     16.559&      20.178&       1.0&     15.602&      22.893&      18.370\\
    & $T_s/T$      &       1.0&      1.380&       2.664&       1.0&      1.948&       2.102&       2.926\\
    \hline\hline
  \end{tabular}
  \caption{Performance of the NMLM solver for the Poiseuille flow with $M=4,5$.}
  \label{tab:poiseuille:nmm-M45}
\end{table}

\begin{table}[!ht]
  \centering\footnotesize
  \begin{tabular}{c|c||ccccccc}
    \hline\hline
    \multicolumn{2}{c||}{} & \multicolumn{6}{c}{$m_{\sss l-1} = m_{\sss l} - 1$}\\
    \hline
    \multicolumn{2}{c||}{$L+1$} & 2 & 3 & 4 & 5 & 6 & 7 & 8 \\
    \hline\hline
    \multirow{4}{*}{\begin{sideways}$N=128$\end{sideways}} 
    & $K$          &         763&         566&         441&         352&         286&         233&         188\\
    & $T$       &    1579.371&    1469.207&     873.198&    1246.409&    1037.293&     872.312&     705.808\\
    & $K_s/K$      &      14.667&      19.772&      25.376&      31.793&      39.129&      48.030&      59.527\\
    & $T_s/T$      &       1.497&       1.610&       2.708&       1.897&       2.280&       2.711&       3.351\\
    \hline\hline
    \multirow{4}{*}{\begin{sideways}$N=256$\end{sideways}} 
    & $K$          &        1680&        1247&         970&         776&         630&         514&         414\\
    & $T$       &    6771.226&    6871.824&    5150.564&    3965.397&    3738.700&    3824.224&    3020.587\\
    & $K_s/K$      &      14.674&      19.769&      25.414&      31.768&      39.130&      47.961&      59.546\\
    & $T_s/T$      &       1.864&       1.837&       2.451&       3.183&       3.376&       3.301&       4.179\\
    \hline\hline
    \multirow{4}{*}{\begin{sideways}$N=512$\end{sideways}} 
    & $K$          &        3560&        2642&        2056&        1646&        1336&        1089&         877\\
    & $T$       &   30046.298&   24663.853&   19606.473&   19143.007&   15356.107&   11740.206&   11755.468\\
    & $K_s/K$      &      14.678&      19.779&      25.416&      31.747&      39.113&      47.984&      59.584\\
    & $T_s/T$      &       1.608&       1.959&       2.465&       2.524&       3.147&       4.116&       4.111\\
    \hline\hline
  \end{tabular}
  \caption{Performance of the NMLM solver for the Poiseuille flow with $M=10$ (part I).}
  \label{tab:poiseuille:nmm-M10-p1}
\end{table}

\begin{table}[!ht]
  \centering\footnotesize
  \begin{tabular}{c|c||cccc|cc|c}
    \hline\hline
    \multicolumn{2}{c||}{} & \multicolumn{4}{c|}{$m_{\sss l-1} = m_{\sss l} - 2$} & \multicolumn{2}{c|}{$m_{\sss l-1} = \lceil m_{\sss l} / 2 \rceil$} & \\
    \hline
    \multicolumn{2}{c||}{$L+1$} & 2 & 3 & 4 & 5 & 2 & 3 & 1\\
    \hline\hline
    \multirow{4}{*}{\begin{sideways}$N=128$\end{sideways}} 
    & $K$          &         722&         498&         346&         217&         569&         340&     11191\\
    & $T$       &    1683.538&    1014.106&     753.820&     483.238&     881.392&     371.769&  2364.837\\
    & $K_s/K$      &      15.500&      22.472&      32.344&      51.571&      19.668&      32.915&       1.0\\
    & $T_s/T$      &       1.405&       2.332&       3.137&       4.894&       2.683&       6.361&       1.0\\
    \hline\hline
    \multirow{4}{*}{\begin{sideways}$N=256$\end{sideways}} 
    & $K$          &        1590&        1098&         761&         478&        1253&         750&     24652\\
    & $T$       &    6355.823&    4845.505&    2014.658&    1982.888&    2369.244&    2227.680& 12623.650\\
    & $K_s/K$      &      15.504&      22.452&      32.394&      51.573&      19.674&      32.869&       1.0\\
    & $T_s/T$      &       1.986&       2.605&       6.266&       6.366&       5.328&       5.667&       1.0\\
    \hline\hline
    \multirow{4}{*}{\begin{sideways}$N=512$\end{sideways}} 
    & $K$          &        3370&        2326&        1613&        1014&        2656&        1589&     52255\\
    & $T$       &   21037.239&   13185.407&   10524.655&    8432.436&   14278.116&    7782.965& 48320.953\\
    & $K_s/K$      &      15.506&      22.466&      32.396&      51.534&      19.674&      32.885&       1.0\\
    & $T_s/T$      &       2.297&       3.665&       4.591&       5.730&       3.384&       6.209&       1.0\\
    \hline\hline
  \end{tabular}
  \caption{Performance of the NMLM solver for the Poiseuille flow with $M=10$ (part II).}
  \label{tab:poiseuille:nmm-M10-p2}
\end{table}

\newcommand\drawPoiseuilleNMMHistory[2]{\begin{figure}[!htb]
  \centering
  {\includegraphics[width=0.5\textwidth]{poiseuille_nmm_res_iters_N#1_M#2.eps}}\hfill
  {\includegraphics[width=0.5\textwidth]{poiseuille_nmm_res_cputime_N#1_M#2.eps}}
  \caption{Convergence history of the NMLM solver for the Poiseuille
    flow with $M=#2$ on a uniform grid of $N=#1$.}
  \label{fig:poiseuille:nmm-res-history-N#1-M#2}
\end{figure}
}
\drawPoiseuilleNMMHistory{512}{4}
\drawPoiseuilleNMMHistory{512}{5}
\drawPoiseuilleNMMHistory{512}{10}

\section{Concluding remarks}
\label{sec:conclusion}

The acceleration for the steady-state computation of the high-order
moment model by using the lower-order model correction has been
investigated in this paper. A nonlinear multi-level moment solver
which has unified framework for the moment model of arbitrary order is
then developed. The convergence rate would be improved as the total
levels of the NMLM solver increases. It is demonstrated by numerical
experiments of two benchmark problems that the proposed NMLM solver
improves the convergence rate significantly and the total
computational cost could be saved a lot, in comparison to the single
level solver. Three order reduction strategies for the lower-order
model correction are also considered. It turns out that the most
efficient strategy is $m_{\sss l-1} = \lceil m_{\sss l} / 2 \rceil$,
the second is $m_{\sss l-1} = m_{\sss l} - 2$, and the third is
$m_{\sss l-1} = m_{\sss l} - 1$.

It should be pointed out that the NMLM solver does not as efficient as
the nonlinear multigrid solver developed in \cite{hu2014nmg}. However,
we have that the spatial grid for our NMLM solver is unchanged at each
level, and the acceleration ratio obtained by the NMLM solver would be
maintained on different spatial grid. Then a natural way of obtaining
a more efficient steady-state solver might be to combine both NMLM
iteration and nonlinear multigrid iteration together. This will be
investigated in our futural work.

\section*{Acknowledgements}
The research of Z. Hu is partially supported by the Natural
  Science Foundation of Jiangsu Province (BK20160784) of China, and
the Hong Kong Research Council ECS grant No.~509213 during his
postdoctoral fellow at the Hong Kong Polytechnic University. The
research of R. Li is supported in part by the National Science
Foundation of China (11325102, 91330205). The research of Z. Qiao is
partially supported by the Hong Kong Research Council ECS grant
No.~509213.


\bibliographystyle{plain}
\bibliography{steadystate,article,tiao}
\end{document}